
\input amstex
\loadbold

\define\aap{Ann.\ Appl.\ Probab.}
\define\ap{Ann.\ Probab.}
\define\mprf{Markov Process.\ Related Fields}
\define\ptrf{Probab.\ Theory Related Fields}

\define\zwvg{Z.\ Wahrsch.\ Verw.\ Gebiete}

\define\spa{Stochastic Process.\ Appl.}

\define\superzd{^{\raise1pt\hbox{$\scriptstyle {\Bbb Z}^d$}}}
\define\xsubvn{x_{\lower1pt\hbox{$\scriptstyle V(n)$}}}
\define\xssubvn{x_{\lower1pt\hbox{$\scriptscriptstyle V(n)$}}}

\define\ssub#1{_{\lower1pt\hbox{$\scriptstyle #1$}}} 
\define\ssubb#1{_{\lower2pt\hbox{$\scriptstyle #1$}}} 
\define\Ssub#1{_{\lower1pt\hbox{$ #1$}}} 
\define\Ssu#1{_{\hbox{$ #1$}}} 
\define\Ssubb#1{_{\lower2pt\hbox{$ #1$}}} 
\define\sssub#1{_{\lower1pt\hbox{$\scriptscriptstyle #1$}}} 
\define\sssubb#1{_{\lower2pt\hbox{$\scriptscriptstyle #1$}}} 
\define\sssu#1{_{\hbox{$\scriptscriptstyle #1$}}}                                              
                                     
\define\<{\langle}
\define\>{\rangle}

\define\({\left(}
\define\){\right)}
\define\[{\left[}                                   
\define\]{\right]}                                            
\define\lbrak{\bigl\{}
\define\rbrak{\bigr\}}
\define\lbrakk{\biggl\{}
\define\rbrakk{\biggr\}}


\define\osup#1{^{\raise1pt\hbox{$ #1 $}}}  
\define\osupm#1{^{\raise1pt\hbox{$\mskip 1mu #1 $}}}  
\define\osupp#1{^{\raise2pt\hbox{$ #1 $}}}  
\define\osuppp#1{^{\raise3pt\hbox{$ #1 $}}}  
\define\sosupp#1{^{\raise2pt\hbox{$\ssize #1 $}}}  
\define\sosup#1{^{\raise1pt\hbox{$\scriptstyle \mskip 1mu #1 $}}} 
\define\ssosup#1{^{\raise1pt\hbox{$\scriptscriptstyle #1 $}}} 
\define\ssosupp#1{^{\raise2pt\hbox{$\scriptscriptstyle #1 $}}} 
\define\ssosuppp#1{^{\raise3pt\hbox{$\scriptscriptstyle #1 $}}}



\define\That{{\widehat T}}
\define\ehat{{\hat e}}
\define\that{{\widehat t}}

\define\ybar{{\overline y}}

\define\Pbar{{\overline P}}

\define\zetatil{{\tilde \zeta}}

\define\sigmatil{{\tilde \sigma}}

\define\OObar{{\overline \Omega}}

\define\pihat{{\widehat \pi}}

\define\oobar{{\overline \omega}}

\define\FFbar{{\overline \FF}}

\define\Rtil{{\widetilde R}}

\define\Ztil{{\widetilde Z}}


\define\sle{\preccurlyeq} 

\define\e{\varepsilon}
\define\a{\alpha}                   
  
\define\ind{{\boldkey 1}} 
\define\inte{\text{\rm int}\,}

\define\epi#1{\text{epi}\,#1}

\define\pikkuhyppy{\vskip .1in}


\define\mmQ{{\bold Q}}

\define\mmN{{\bold N}}

\define\mmR{{\bold R}}

\define\mmZ{{\bold Z}}


\define\SS{{\Cal S}}
\define\cAA{{\Cal A}}    
\define\BB{{\Cal B}}
\define\CC{{\Cal C}}

\define\FF{{\Cal F}}
\define\GG{{\Cal G}}

\define\II{{\Cal I}}

\define\NN{{\Cal N}}

\define\TT{{\Cal T}}      
\define\WW{{\Cal W}}      
      
\define\ZZ{{\Cal Z}}  

\define\BBtil{{\widetilde \BB}}

\define\TTtil{{\widetilde \TT}}   
\define\WWtil{{\widetilde \WW}}   



\documentstyle{amsppt}
\magnification=\magstep1
\document     
\baselineskip=12pt
\pageheight{43pc} 

\NoBlackBoxes

\centerline{\bf Strong Law of Large Numbers }

\pikkuhyppy

\centerline{\bf  for the Interface in Ballistic Deposition }

\pikkuhyppy

$$\text{1999}$$

\hbox{}

\centerline{  Timo  Sepp\"al\"ainen
\footnote""{ Research partially supported by NSF grant DMS-9801085. }}
\hbox{}
\centerline{Department of Mathematics}
\centerline{Iowa State University}
\centerline{Ames, Iowa 50011, USA}
\centerline{seppalai\@iastate.edu}
\vfil

\flushpar
{\it Summary.} We prove a hydrodynamic limit for 
ballistic deposition on a multidimensional lattice.
In this growth model particles rain down at random and 
stick to the growing cluster at the first point of contact. 
The theorem is that 
if the initial random interface converges to a 
deterministic macroscopic function, then 
at later times the height of
the scaled interface converges to the viscosity
solution of a Hamilton-Jacobi equation. 
The proof idea is to decompose 
 the interface into  the shapes 
that grow from individual seeds
of the initial interface. This decomposition 
converges to a variational formula that
defines viscosity solutions of the macrosopic equation. 
The technical side of the proof involves subadditive 
methods and large deviation bounds for related 
first-passage percolation processes. 


\hfil

\hfil


\hfil

\flushpar

\hfil

\flushpar
Short Title: Ballistic deposition

\break

\head 1.  Introduction \endhead

In a ballistic deposition model particles fall on 
a surface, find a location where they attach, and
become part of the growing cluster. Depending on 
the rules chosen, the particle may stick to the 
first point of contact, 
or it may ``roll downhill'' and attach itself to the 
first stable location it encounters. We study the 
version where  particles stick to the first 
point of contact. A ballistic deposition
 model is {\it flux limited}
as opposed to {\it reaction limited} because 
the rate of growth is limited by the availability of
material rather than by the availability of growth sites. 

Another
qualitative distinction between growth models is 
the division in {\it local} and {\it nonlocal} models. 
An example of a local model is one where
the rate of attaching new particles depends only 
on the states of some 
finite number of neighboring sites. A highly
nonlocal model is diffusion limited
aggregation where
 the chances of a particle attaching at any particular
site are influenced by 
 far-away parts of the growing cluster.
Ballistic deposition has a nonlocal aspect because
it can build overhangs that extend far sideways
to shade parts of the interface that then no longer 
receive the  particle flux. However, the {\it height}
of the interface has local dynamics, and it is this
height process that we study in our paper.

 The shadowing effect and low
atomic mobility (stick to the first point of 
contact) may play a role in actual
deposition processes, so there is 
 physical motivation for these rules (see \cite{16}). 
We refer the reader to the survey article of Krug and Spohn
\cite{18} 
for a general discussion of growth models and for 
references to the physics literature. 

We prove a {\it hydrodynamic scaling limit} for 
a ballistic deposition  process on an infinite 
$d$-dimensional cubic lattice. 
Particle deposition happens independently 
at all sites, governed by exponential waiting times. 
Once deposited, particles never leave the surface,
so in particle system jargon this process
is {\it totally asymmetric}. The 
correct hydrodynamic scaling shrinks space
and speeds up time by the same factor $n$. The 
theorem is a law of large numbers: as $n\to\infty$,
the height of the randomly evolving surface converges to the 
solution of a first-order partial differential
equation of the Hamilton-Jacobi type. For general
treatments of hydrodynamic limits,
the reader may consult references \cite{5}, \cite{15} and  
\cite{29}. 

From a statistical mechanics point of view, 
our result is a rigorous microscopic derivation 
of the expected macroscopic theory. The macroscopic
growth velocity of the height function
 $\psi$ is determined
by the local slope: $\partial\psi/\partial t=f(\nabla\psi)$.
The function $f$ that gives the dependence
is the Legendre conjugate of the stationary shape $g$
of a cluster grown from  a seed:
$f(p)=\sup_x\{x\cdot p+g(x)\}$. As explained in 
\cite{18}, this can be viewed as a 
 Wulff construction for a growing shape. 

The essence of the proof is to 
construct the process so that 
 a supremum
of ballistic deposition processes is again a ballistic
deposition process. The proof works equally well 
for discrete time in which case the waiting times
of deposited particles are geometric rather than 
exponential. In fact the Markovian nature of the 
dynamics is not needed at all. The proof can be carried out for 
arbitrary waiting times, but we have not bothered with such 
generality  here. An example
can be found in \cite{26}. 

The paper is organized as follows. Section 2 defines
the model more precisely and states the theorems.
We discuss briefly viscosity solutions of
Hamilton-Jacobi equations, and suggest
 open problems from the physics
literature that may be amenable to rigorous
progress. In Section 3 we do the standard graphical
construction (see for example \cite{6}, \cite{9}, \cite{11}, or 
\cite{19}) to show that the process can be rigorously
defined on a probability space that supports the 
Poisson jump time processes. And we construct the 
coupling that expresses the process as the
supremum of the shapes from individual seeds of the 
initial interface. 

In Section 4 we prove exponential 
large deviation estimates for a first-passage site 
percolation process. We need these because the
lateral (sideways) growth of ballistic deposition from
a seed is equivalent to a first-passage percolation 
problem.  The bounds we need are of a standard type, 
but they are
not available in the literature in the right form. 
 We use a result of Talagrand \cite{30}
and some ideas from Kesten \cite{14} to derive the
inequalities.
Grimmett and Kesten \cite{10} have proved 
such bounds for growth along coordinate axes, but 
we need to control growth in all space directions. 
Generalizing the block argument of \cite{10} appeared 
harder than applying the very general tools of \cite{30}.

Theorem 1 is proved in Section 5, and
in Section 6 a technical extension of Theorem 1 is
proved where the initial seed and initial time are
translated as the limit is taken. This is required
for the proof of Theorem 2, which is the content 
of Section 7. 

\subsubhead Notational conventions 
\endsubsubhead
Frequently used notation is summarized here for the 
reader's convenience. $C$, $C_1$, $C_2$ are
constants whose exact values are immaterial and whose
values may change from one inequality to the next.
 $\mmR_+=\{x\in\mmR:x\ge 0\}$ is the 
set of nonnegative real numbers, and similarly for $\mmQ_+$ and 
$\mmZ_+$. $\mmN=\{1,2,3,\ldots\}$ is the set of natural numbers. 
{\it Sites} or points of
$\mmZ^d$ are denoted by $u,v,w,z$.
$0$ is the origin of $\mmZ^d$, and $\NN$ is the set of nearest
neighbors of the origin in $\mmZ^d$, in other words, the $2d$ 
sites at $\ell^1$-distance 1 from the origin. To
distinguish points of 
$\mmZ^{d+1}$ from those of $\mmZ^d$,
we  call points of $\mmZ^{d+1}$ {\it cells} and denote
them by $(u,k)$ where $u\in\mmZ^d$ and $k\in\mmZ$.  
In particular,
$(0,0)$ is the origin of $\mmZ^d\times\mmZ_+$.
 For a real number $x$, $[x]$ is the maximal
integer $n$ subject to $n\le x$. For $x=(x_1,\ldots,x_d)\in\mmR^d$ the
site $[x]=([x_1],\ldots,[x_d])$, and the $\ell^1$ norm
is $|x|=|x_1|+\cdots+|x_d|$. For $x,y\in\mmR$,
 $x\vee y=\max\{x,y\}$,
 $x\wedge y=\min\{x,y\}$. 
$\|g\|_\infty=\sup_x|g(x)|$ for any function $g$. 

 $\ind_A$ and $\ind\{A\}$ denote the indicator random variable
of the event $A$. 
 An Exp($\beta$) 
 random
variable $X$ satisfies $P(X>t)=\exp(-\beta t)$ for $t\ge 0$.
$S^1_n$ stands for a sum of $n$ i.i.d.\ Exp(1) random 
variables. The standard notion of stochastic dominance is expressed
by $X\sle Y$ which is equivalent to $P(X\le t)\ge P(Y\le t)$
for all $t$. $\theta_s$ is a time translation on Poisson
point processes that translates points $r$ to $r-s$. 
In other words, reading $\theta_s\omega$ from time $0$
onwards is the same as  reading $\omega$ from time $s$
onwards. $B_0\subseteq\mmR^d$ is the convex compact limiting shape for
first-passage site percolation on $\mmZ^d$  with Exp(1) waiting times.

\hbox{}

\head 2. The results \endhead

We start with a description of the 
 ballistic deposition process. Fix a positive
integer $d$. 
Imagine first an arbitrary subset of $\mmZ^{d+1}$ 
of occupied {\it cells}. Cells are simply  points of $\mmZ^{d+1}$
and  denoted by $(u,h)$, where $u\in\mmZ^d$ is a {\it site}
and $h\in\mmZ$ is a {\it height}. The {\it cluster} 
(subset of occupied cells) grows through 
deposition events: over each site $u\in\mmZ^d$,  
particles rain down randomly 
at exponential rate 1, independently of all other sites. 
A particle descending down  sticks to the first spot where it
touches the existing cluster. We imagine that
the particles are exactly the size of a unit
cube in $d+1$ dimensions. Over a site $u$ a particle
instantaneously drops down from height $h=\infty$
to the highest cell adjacent to the 
existing cluster, and then this cell becomes 
occupied and joins the cluster. 
(Cells are {\it adjacent} if their $\ell^1$ distance
is 1.)
As the process evolves, a porous structure is 
generated that grows upward in $d+1$ dimensions. 
 
We follow the evolution of the top surface
of the cluster. For each site $u$, the {\it height}
 $\sigma_u$ is 
the maximal $h$ such that cell $(u,h)$ is occupied. 
We permit the values  $\sigma_u=\pm\infty$.
$\sigma_u=-\infty$ means that no cell in the 
column $\{(u,h):h\in\mmZ\}$ is occupied. 
The state of the process is the configuration 
$\sigma=(\sigma_u:u\in\mmZ^d)$ of height 
variables, and  the state space is
  $\SS=(\mmZ\cup\{\pm\infty\})^{\mmZ^d}$. 
We write $\sigma(t)=(\sigma_u(t))$ 
for the height process, where $t\ge 0$
denotes time. 

 The time evolution
of the interface is determined by a collection
$\TT=\{\TT^u:u\in\mmZ^d\}$ of rate 1 Poisson point
processes on the time line $(0,\infty)$. At each epoch 
of $\TT^u$  a  particle is deposited 
 above site $u$. In terms of the height
variables $\sigma_u(t)$ a deposition event has a simple
 expression: Let $\NN$ denote the set of $2d$  nearest
neighbors of the origin in $\mmZ^d$.
 If $r$ is an epoch in Poisson
process $\TT^u$, the height above site $u$ jumps
at time $r$ according to this formula:
$$\sigma_u(r)= \max\bigl[ \sigma_u(r-)+1\,,\,
\max\{ \sigma_{u+z}(r-): z\in\NN\}\bigr]\,.
\tag 2.1
$$
The maximum in (2.1) 
has the effect  that the deposited particle
sticks to the highest unoccupied cell above site $u$
that is adjacent to the existing cluster. 
If no cell above site $u$ is adjacent to the cluster
at time $r-$, the right-hand side of (2.1) equals $-\infty$.
This means that the deposited particle is lost
because it cannot stick to the cluster.
In Section 3 we construct this process rigorously
starting from an arbitrary initial interface $\sigma(0)$. 

The simplest  
ballistic deposition process starts from
a single occupied cell, or {\it seed}. The prototypical
one starts from a seed at the origin. For this process
we reserve the symbol $Z$ and use the symbol $\sigma$ 
for the general ballistic deposition process. 
So at time 0,
$Z_0(0)=0$, and $Z_u(0)=-\infty$ for all sites $u\ne 0$. 
Otherwise $Z(\cdot)$ evolves as specified above. 

Let $\BB(t)$ denote the set of sites above which $Z$
has an occupied cell  by time $t$:
$$\BB(t)=\{u\in\mmZ^d: Z_u(t)\ge 0\}.
\tag 2.2 
$$
Since the original seed is at the origin, $Z_u(t)\ge 0$
is equivalent to $Z_u(t)>-\infty$. 
The cluster $\BB(t)$ is the view of the process $Z(t)$
from above, by projecting $\mmZ^{d+1}$ 
onto $\mmZ^d$ by the map $(u,h)\mapsto u$. 
The ballistic deposition rules imply that $\BB(t)$
grows according to this rule: Each site adjacent 
to the existing cluster joins independently at rate 1. 
[Because particles rain down  at rate 1 above each site $u$,
and they have a chance of sticking iff some cell above
an adjacent site is already occupied.]
Thus $\BB(t)$ is a 
familiar growth model, namely {\it first-passage 
percolation}, or a version of the Eden growth model.
It  satisfies a law of large 
numbers. There is a closed, convex 
deterministic set $B_0\subseteq\mmR^d$ 
with 
nonempty interior such that this holds almost 
surely: for any $\e>0$,
$$\bigl[ t(1-\e)B_0\bigr]\cap \mmZ^d\subseteq\BB(t)\subseteq
\bigl[ t(1+\e)B_0\bigr]\cap \mmZ^d 
\tag 2.3
$$
for all large enough $t$. 
In first-passage percolation literature 
 the random waiting times are
usually attached to the edges, but for our $\BB(t)$ the
waiting times
are attached to the sites. In Section 4 we prove 
some large deviation estimates for $\BB(t)$. 

Over the set $B_0$ the scaled
 ballistic deposition process $t^{-1}Z(t)$
 approaches a limiting shape.
Let $\inte B_0$ denote the topological interior of 
 $B_0$. 
 We prove in Section 5 this
theorem: 

\proclaim{Theorem 1} There exists a 
bounded positive function
$g$ defined on the open convex set $\inte B_0$ such that
this law of large numbers holds: Outside an 
event of probability zero, 
$$\lim_{n\to\infty} \frac1n Z_{[nx]}(nt)=tg(x/t) 
\tag 2.4 
$$
for all $x\in\mmR^d$ and $t>0$ such that $x/t\in\inte B_0$. 
Furthermore, $g$ is continuous, concave, and invariant under
 permutations of the coordinate axes and reflections
about the origin. 
\endproclaim

The main theorem is a scaling limit from a general 
initial interface. Suppose we have  a sequence $\sigma^n$
of ballistic deposition processes. The initial
macroscopic interface is given by a function $\psi_0$ 
on $\mmR^d$. We consider three different sets of
hypotheses.

\hbox{}

\subhead Assumption A \endsubhead 
$\psi_0$ is a continuous $[-\infty,+\infty]$-valued
 function on $\mmR^d$. 
The ballistic
deposition processes $\sigma^n(t)$ are constructed
on a common probability space, and all processes 
 use the same version of the Poisson jump time
processes. There is a countable dense subset
$Y_0\subseteq\mmR^d$ such that 
$$\lim_{n\to\infty} \frac1n \sigma^n_{[ny]}(0)= \psi_0(y)
\tag 2.5
$$
 almost surely for each $y\in Y_0$, and 
each  $y\in Y_0$ has arbitrarily small closed
neighborhoods $V$ such that almost surely 
$$\limsup_{n\to\infty} \frac1n \cdot\max_{u\in nV}\sigma^n_{u}(0)\le 
\sup_{y\in V}\psi_0(y)\,.
\tag 2.6
$$

\hbox{}

\subhead Assumption B \endsubhead 
$\psi_0$ is an arbitrary $[-\infty,+\infty]$-valued
 function on $\mmR^d$. 
The ballistic
deposition processes $\sigma^n(t)$ are constructed
on a common probability space, and all processes 
 use the same version of the Poisson jump time
processes. For any fixed $y\in\mmR^d$, the limit  (2.5) holds
almost surely, and there are arbitrarily
small  closed
neighborhoods $V$ of $y$ such that  (2.6) holds almost surely.

\hbox{}

\subhead Assumption C \endsubhead 
Again $\psi_0$ is an arbitrary $[-\infty,+\infty]$-valued
 function on $\mmR^d$.   
 (2.5) and (2.6) hold
in probability in this sense: Let
$(\Omega_n,\FF_n,P_n)$ be the probability space of the
$n$th process $\sigma^n$. For any $y\in\mmR^d$ and  $\e>0$,
$$\lim_{n\to\infty}P_n\biggl(\, \bigl|
n^{-1} \sigma^n_{[ny]}(0)- \psi_0(y) \bigr|\ge \e\biggr)=0\,,
\tag 2.7
$$
 and  there are arbitrarily small closed
neighborhoods $V$ of $y$ such that 
$$\lim_{n\to\infty} P_n\biggl(\, 
\max_{u\in nV}\sigma^n_{u}(0)\ge 
n\cdot\sup_{y\in V}\psi_0(y)+n\e \biggr)=0\,.
\tag 2.8
$$

\hbox{}

For $x\in\mmR^d$ and $t>0$, define 
$$\psi(x,t)=\sup_{y\in x+t\cdot\inte B_0}\lbrakk \psi_0(y)
+t\,g\biggl(\frac{x-y}t\biggr)\rbrakk\,,
\tag 2.9
$$
and $\psi(x,0)=\psi_0(x)$. 

\proclaim{Theorem 2} 
{\rm (a)} Strong Laws of Large Numbers: 
Under Assumption A we have this convergence, 
simultaneously for all
 $x\in\mmR^d$ and $t>0$, outside 
a single exceptional event of probability zero: 
$$\lim_{n\to\infty}\frac1n\sigma_{[nx]}^n(nt)=\psi(x,t)\,.
\tag 2.10
$$
Under Assumption B  the 
limit in {\rm (2.10)} holds almost surely for all $(x,t)$ at which 
$\psi$ is $t$-continuous from the right, $\psi(x,t)=\psi(x,t+)$. 

{\rm (b)} Weak Law of Large Numbers: Under Assumption C the 
limit in {\rm (2.10)} holds in probability for all
 $(x,t)$ at which 
$\psi$ is $t$-continuous from the right.  In other words, 
for all such  $(x,t)$, 
$$\lim_{n\to\infty}P_n\biggl(\, \bigl|
n^{-1} \sigma^n_{[nx]}(nt)- \psi(x,t) \bigr|\ge \e\biggr)=0
\tag 2.11
$$
for all $\e>0$. 
\endproclaim

\hbox{}

\flushpar
{\bf Remarks and extensions}

\hbox{}

\flushpar
{\it 2.1 The assumptions.}
 Some uniformity assumption such as (2.6)
 is needed for the result. Consider this
example in dimension $d=1$: The initial
interface is 
$$\sigma^n_u(0)=\cases n, &u=1\\
           0, &u\ne 1\,.
\endcases
$$
The limit (2.5) is satisfied with $\psi_0(x)\equiv 0$,
and then from (2.9) $\psi(x,t)=tg(0)$. 
But the spike in $\sigma^n_u(0)$ at site $u=1$ implies that 
$$n^{-1}\sigma_0(nt)\ge 1+n^{-1}Z^1_{0}(nt)
\longrightarrow 1+tg(0)\qquad\text{as $n\to\infty$,}
$$
where we wrote $Z^1$ for the process that starts from a seed
in cell $(1,0)$. 

The point of separating Assumptions A and B is that
continuity of $\psi_0$ guarantees continuity of $\psi(x,t)$
[Lemma 7.1 in Section 7]. 
Then  the limit (2.10) can be proved almost surely,
 {\it simultaneously for all} $(x,t)$, even though the
hypothesis requires the limit only on a countable dense
 set $Y_0$. The limit under Assumption
B generalizes Theorem 1. The right-continuity 
in $t$ at $(x,t)$, namely $\psi(x,t)=\psi(x,t+)$, follows
from $x$-continuity at $(x,t)$, as can be verified from (2.15)
below. 

The nonlocal shadowing effect restricts the 
law of large numbers to continuity points of the 
macroscopic surface. Consider the process $Z$ growing
from a seed at $(0,0)$ in dimension $d=1$. The limit (2.4)
cannot possibly hold at $x=t$ because whether site
$[nt]$ is shadowed by the top surface varies with
the fluctuations of the right edge of the growing 
cluster: $Z_{[nt]}(nt)=-\infty$ or $\ge 0$ 
depending on whether $S^1_{[nt]}>nt$ or $\le nt$,
where $S^1_{[nt]}$ is a sum of $[nt]$ i.i.d.\ Exp(1) random 
variables and represents the time it takes for the
cluster to grow $[nt]$ lattice units to the right. 

The proof of Theorem 2 depends on two essential things:
 a coupling (Lemma 3.3) that expresses the general 
process as a supremum of processes of the $Z$ type, and 
the limit in Theorem 1 for the $Z$ process. The coupling
is basically a combinatorial property of 
the paths. As long as  pathologies such as 
jump times accumulating at a point are ruled out,
the coupling  is independent of the probability
distribution on the paths. The limit in Theorem 1
on the other hand comes from an ergodic theorem 
and is closely tied to the probabilities
of the evolution of $Z$. Thus our proof of Theorem 2 
works for any reasonable  jump rates under 
which Theorem 1 can be proved. For example, we could 
take {\it quenched} random jump rates where deposition
at site $u$ happens at rate $\alpha_u$, and the 
rates $\{\alpha_u\}$ are i.i.d.\ random variables 
(for examples in other contexts, 
see \cite{2}, \cite{26} and \cite{27}). Another
possibility  would be to let a sufficiently regular
function $\alpha(x)$ on 
$\mmR^d$ determine the rates so that deposition  at site $u$ for 
process $\sigma^n$ 
 happens at rate $\alpha(u/n)$ (for examples, see
\cite{1} and \cite{21}).

\hbox{}

\flushpar
{\it 2.2 Viscosity solutions of Hamilton-Jacobi equations.}
Consider  the Hamilton-Jacobi equation
$$\align
\psi_t-f(\nabla\psi)&=0 \tag 2.12\\
\psi(x,0)&=\psi_0(x)\tag 2.13
\endalign
$$
where $f$ is some continuous function on $\mmR^d$. 
A function $\psi(x,t)$ on $\mmR^d\times\mmR_+$ that 
satisfies the initial condition (2.13) is  a
{\it viscosity solution} of (2.12)--(2.13) if the following
holds for all continuously differentiable functions
$\phi$ on $\mmR^d\times(0,\infty)$: 
if $\psi-\phi$ has a local maximum at $(x_0,t_0)$, then
$$\phi_t(x_0,t_0)-f\bigl(\nabla\phi(x_0,t_0)\bigr)\le 0\,,$$
and if $\psi-\phi$ has a local minimum at $(x_0,t_0)$, then
$$\phi_t(x_0,t_0)-f\bigl(\nabla\phi(x_0,t_0)\bigr)\ge 0\,.$$
The notion of viscosity solution is due to Crandall and
Lions \cite{4}. Properties of viscosity solutions  
of Hamilton-Jacobi equations are developed 
in  Crandall, Evans and Lions
\cite{3}, and in the textbook of Evans \cite{8}.  

Equation (2.9) is known as a {\it Hopf-Lax
formula} in the partial differential equations literature.
Let  $f$ be the (negative of the) Legendre conjugate of $g$:
$$f(u)=\sup_{x\in\inte B_0}\{ u\cdot x+g(x)\}\,,\qquad u\in\mmR^d\,.
\tag 2.14
$$
Since $B_0$ is a compact set and $g$ is bounded, one can check
that $f$ is finite and Lipschitz continuous on all of $\mmR^d$. 
Then an  argument in \cite{8} (proof of Theorem 3 in Section 10.3.4) 
 shows
that $\psi(x,t)$ defined by (2.9) is a viscosity solution 
of (2.12)--(2.13). 
By applying results from the p.d.e.\ literature, we can 
refine Theorem 1 with a 
 uniqueness statement: 

\proclaim{Theorem 3}
Suppose Assumption A is in force, 
and that additionally the initial macroscopic 
profile $\psi_0$ is uniformly 
continuous on $\mmR^d$. Fix a finite time horizon
$T<\infty$. As in Theorem 1 the 
strong law of large numbers  {\rm (2.10)} is valid.
On $\mmR^d\times[0,T]$  the limit $\psi(x,t)$ 
 is the 
unique uniformly continuous viscosity solution 
  of the Hamilton-Jacobi equation {\rm (2.12)--(2.13)}
whose Hamiltonian $-f$ is defined by {\rm (2.14)}. 
\endproclaim

\demo{Proof of Theorem 3}  
The point is that the additional assumption
of uniform continuity on $\psi_0$ permits us to
make the  uniqueness assertion about $\psi$. Uniqueness
theorems for unbounded viscosity solutions have been
proved by  Ishii \cite{12}, whose
   Theorem 2.1 
states that equation (2.12)--(2.13)  with continuous $f$
has a unique uniformly continuous viscosity solution 
 on $\mmR^d\times[0,T]$.   Assume Theorem 2. 
To prove Theorem 3, 
 we therefore need to check the uniform continuity
of $\psi $ defined by (2.9), assuming that $\psi_0$ 
 is uniformly continuous. Here is an outline of the argument. 

First check that formula (2.9) operates like a semigroup:
Once $\psi$ has been defined by (2.9), it follows for all
$0<s<t$ that
$$\psi(x,t)=\sup_{y\in x+(t-s)\cdot\inte B_0}
\lbrakk \psi(y,s)
+(t-s)g\biggl(\frac{x-y}{t-s}\biggr)\rbrakk\,.
\tag 2.15
$$
This bounds the growth of $\psi(x,t)$ 
in time. Let $b_0=\sup\{|y|:y\in B_0\}$. Fix $x$ and 
$s<t$. (2.15) gives
$$\aligned
0&\le \psi(x,t)-\psi(x,s)\\
&\le \sup_{|y'-y''|\le b_0(t-s)}\bigl|\psi(y',s)-\psi(y'',s)\bigr|
+(t-s)\|g\|_\infty\,.
\endaligned
\tag 2.16
$$
On the other hand, for a fixed time $t$ the definition (2.9)
directly bounds the variation of $\psi(\cdot\,, t)$:
For any $x_1,x_2\in\mmR^2$ and $t>0$, 
$$\bigl|\psi(x_1,t)-\psi(x_2,t)\bigr|
\le  \sup_{|y'-y''|\le |x_1-x_2|}\bigl|\psi_0(y')-\psi_0(y'')\bigr|\,.
\tag 2.17
$$
Now the uniform continuity of $\psi$ on all of
$\mmR^d\times\mmR_+$ follows: 
$$\aligned
&\bigl|\psi(x_1,t_1)-\psi(x_2,t_2)\bigr|\\
&\le \bigl|\psi(x_1,t_1)-\psi(x_2,t_1)\bigr|
+  \bigl|\psi(x_2,t_1)-\psi(x_2,t_2)\bigr|\\
&\le \sup_{|y'-y''|\le |x_1-x_2|}\bigl|\psi_0(y')-\psi_0(y'')\bigr|\\
 &\qquad +  \sup_{|x'-x''|\le b_0(t_2-t_1)}
\bigl|\psi(x',t_1)-\psi(x'',t_1)\bigr|
+(t_2-t_1)\|g\|_\infty\\
&\qquad\qquad\qquad\qquad[\text{ without loss of generality,
assume that $t_1<t_2$ }]\\
&\le 2\cdot \sup\lbrak
\bigl|\psi_0(y')-\psi_0(y'')\bigr|: |y'-y''|\le |x_1-x_2|\vee b_0|t_2-t_1|
\rbrak \\
&\qquad+|t_2-t_1|\cdot\|g\|_\infty\,.
\endaligned
$$
This proves Theorem 3. 
\qed
\enddemo

\hbox{}

\flushpar
{\it 2.3 Statistical mechanics.}
From the point of view of statistical mechanics,
 our paper provides a  rigorous derivation of the 
macroscopic theory  that is  taken as basic in 
the physics literature. According to  this theory, 
macroscopically the interface moves under an
 inclination-dependent
growth  velocity $f(u)$, and $f$ is the Legendre conjugate
of  the cluster shape $g(x)$
that grows from a seed. From this basis 
the physics literature seeks to describe finer properties
of the deposition process. The reader is 
referred to the survey article \cite{18},
and to articles \cite{16}, \cite{17}, and \cite{20}. 

Here we comment on some properties
of the macroscopic objects, and mention open 
problems suggested by the physics
papers. In general, describing 
 $f$ and $g$ is as hard as 
first-passage percolation, since knowing $g$
would imply knowing the first-passage percolation 
shape $B_0$. In dimension $d=1$ the percolation 
question becomes trivial, $B_0=[-1,1]$, so one may
hope to say something about $f$ and $g$ for $d=1$. 
The case $d=1$ is also the focus of the physics
literature. 

 It follows from (2.14) that the 
velocity $f(u)$ is convex and  even [$f(u)=f(-u)$]. 
Consequently it has  a minimum at $u=0$, $f(0)=g(0)$. 
Whether this minimum is strict as expected 
(p.\ 903 in \cite{17}) is a harder question because
that would require showing that $g$ does not have a 
corner at $x=0$. These types of questions are open 
for other interface  models too, except in those 
rare cases where invariant distributions can be 
used to explicitly compute limiting shapes.
For examples, see \cite{23}--\cite{27}. In $d=1$
(2.14) gives linear asymptotics 
$f(u)=|u|+g(1)+o(|u|)$ for large slopes $u\to\pm\infty$. 
The precise nature of the error $o(|u|)$ would be 
of interest. Equivalently, one wants to know the 
asymptotics of $g(x)-g(1-)$ as $x\nearrow 1$. 

The velocity must increase 
with dimension, since  higher 
dimension  means more neighbors to
speed up the growth over a particular site. This is 
easy to check by coupling the $Z$-processes for
 dimensions $d$ and $d+1$ so that the $d$-dimensional
$Z$-process grows over the hyperplane $\{x_{d+1}=0\}$ 
in $\mmZ^{d+1}$. Without any estimates, this gives
$g_{d+1}(x',0)\ge g_d(x')$ for any $x'\in\mmR^d$. 
Consequently, for any $u=(u',u_{d+1})\in\mmR^{d+1}$,
$$\aligned
f_{d+1}(u)&\ge \sup_{x'\in\mmR^d}\{ u'\cdot x'+g_{d+1}(x',0)\}\\
&\ge \sup_{x'\in\mmR^d}\{ u'\cdot x'+g_{d}(x')\}\\
&=f_d(u')\,.
\endaligned
$$
Mean-field analysis in \cite{17} suggests that
$f_d(0)$  grows like $d/\log d$, and 
simulations appear to show a slow convergence toward
mean-field values  as $d\to\infty$. Growth at rate
$d/\log d$ has been verified for first-passage percolation
\cite{13}, so these questions can be investigated
rigorously. 

In a {\it synchronously updated} ballistic deposition
process 
time is discrete, and a rate parameter $p\in(0,1)$ 
is fixed. At each time step $t=1,2,3,\ldots$ , 
 an independent random choice is made at
each site:
with probability $p$ 
the height  is updated according to equation (2.1),
and  with probability
$1-p$ the height remains the same. 
The results of our paper apply to the synchronous
process without changes. The only difference is that
the Poisson point processes $\{\TT^w: w\in\mmZ^d\}$ 
 of jump times are replaced by
Bernoulli processes on the  discrete time line $\mmN=\{1,2,3,\ldots\}$.
In these processes  an event
arrives  at each time with probability $p$, independently of 
everything else. Analogously with the flat edge result 
of Durrett and Liggett \cite{7} for first-passage percolation,
a {\it faceting transition} happens in ballistic deposition for
 large enough $p$. Interestingly, simulations  in \cite{17}
suggest that the velocity $f(u)$ is linear for all
slopes $u\ge 1$
if $p$ is large enough. 

\hbox{}

\head 3. Construction and coupling \endhead

To construct the ballistic deposition process, 
start by giving each site $u\in\mmZ^d$ an independent 
rate 1 Poisson point process $\TT^u$ on the time
line $(0,\infty)$. 
 Fix an initial configuration
$\sigma(0)=(\sigma_u(0):u\in\mmZ^d)\in\SS
\equiv (\mmZ\cup\{\pm\infty\})^{\mmZ^d}$. Informally speaking,
 the construction 
of the dynamics goes as follows: If $r$ is an
epoch (in other words, a point) 
of $\TT^u$, then at time $r$ height variable
 $\sigma_u$ jumps:  
$$\sigma_u(r)= \max\bigl[ \sigma_u(r-)+1\,,\,
\max\{ \sigma_{u+z}(r-): z\in\NN\}\bigr]\,.
\tag 3.1
$$
Recall that  $\NN$ is the set of nearest
neighbors of the origin in $\mmZ^d$. 

To make the 
construction rigorous, we show that there exists a fixed
time $t_0>0$ and a set of Poisson 
processes $\{\TT^u\}$ of full probability
 such that, starting with an arbitrary 
$\sigma(0)\in\SS$, the evolution $\sigma(t)$ can be 
computed for $t\in[0,t_0]$.
 Since $t_0$ is independent
of $\sigma(0)$, the construction can be repeated, starting with
$\sigma(t_0)$, and extended to time interval $[0,2t_0]$. And so on, to
arbitrarily large times. 

Given a fixed number $t_0>0$
and the Poisson processes $\{\TT^u\}$, construct the following
random graph with vertex set $\mmZ^d$: Connect
nearest neighbors  $u$ and $v$
with an edge if either $\TT^u$ or
$\TT^v$ has an epoch in $[0,t_0]$. 

\proclaim{Lemma 3.1} For small enough fixed $t_0>0$, this random
graph has no infinite connected components 
for almost every realization of $\{\TT^u\}$.
\endproclaim

Before proving the lemma, let us see how
the construction problem is solved.
Make these further assumptions on  
the Poisson point processes,  valid for almost every
realization:
\roster
\item "(3.2a)" The $\{\TT^u\}$ are
such that there are no simultaneous jump attempts. 
\item "(3.2b)" Each
$\TT^u$ has finitely may epochs in each bounded time
interval. 
\endroster
 All sites $w$ that can influence the evolution at site $u$ up to time 
$t_0$ are connected to $u$ in the random graph. Since 
$u$ lies in a finite connected component $\CC$, the 
point process $\cup_{w\in\CC}\TT^w$ has only finitely
many epochs in $[0,t_0]$. Consequently the evolution
$\sigma_w(t)$ can be computed for $w\in\CC$ and $t\in[0,t_0]$
from  rule (3.1), by considering 
the finitely many epochs in their temporal order. This
procedure is repeated for all connected components. 

\demo{Proof of Lemma 3.1} By translation invariance, it 
suffices to show that the origin is almost surely
connected to only finitely many vertices.
 $\{u_0,u_1,\ldots,u_n\}$ is a 
self-avoiding path of length $n$  
 in the random graph
 if $u_i\ne u_j$ for $i\ne j$ and if
there is an edge between $u_i$ and $u_{i+1}$ for each $i$.
If the origin is connected to a site $u$ with $|u|\ge L$,
there is a self-avoiding path of length $\ge$ $L$ 
 starting at the origin. The probability
that a  self-avoiding path  of length $2n-1$ 
starts at the origin is at most
$$(2d)^{2n-1}\bigl(1-e^{-2t_0}\bigr)^n\,.
$$
The factor $(2d)^{2n-1}$ is an upper
bound on the number of such paths. If 
$0=u_0,u_1,\ldots,u_{2n-1}$ is such  a path, the $n$
edges $(u_0,u_1)$,  $(u_2,u_3)$, $\ldots$,  $(u_{2n-2},u_{2n-1})$ 
are present independently of each other, and each with
probability $1-e^{-2t_0}$ [at least one
of $\TT^{u_{2i}}$ and  $\TT^{u_{2i+1}}$ must have an epoch in
$[0,t_0]$, and each $\TT^{u_{j}}$ has rate $1$]. 
Pick $t_0$ small enough so that 
$(2d)^2\bigl(1-e^{-2t_0}\bigr)<1$. Then by Borel-Cantelli,
self-avoiding paths from the origin have a finite
upper bound on their length, almost surely. 
\qed
\enddemo

This approach to the construction of a particle system
is due to Harris \cite{11}. Our presentation followed 
\cite{6}. 

Let $(\Omega,\FF,P)$ denote the probability space 
whose sample point $\omega$ represents a realization
of the Poisson processes $\TT=\{\TT^u\}$. 
We constructed  the random path 
$\sigma(\cdot)=(\sigma_u(t):u\in\mmZ^d, t\ge 0)$ 
 as a function of
the given initial state $\sigma(0)$ 
and a sample point $\omega$.  Since the Poisson
processes are Markovian, the process $\sigma(\cdot)$
is a time-homogeneous Markov process.
When the initial interface $\sigma(0)$ is random, 
the underlying probability space is constructed so that
 $\sigma(0)$ and $\TT$ are independent. 

 Formula 
(3.1) preserves ordering, so we get the following
monotonicity lemma, whose proof is left to the
reader:

\proclaim{Lemma 3.2} Suppose 
$\sigma$ and $\rho$ are ballistic
deposition processes 
 constructed on a common probability space so that
they use the same version $\{\TT^u\}$ of the Poisson processes.
 Assume that at time $0$, 
$\sigma_u(0)\ge \rho_u(0)$ for all
$u\in\mmZ^d$, almost surely.
Then almost surely 
$\sigma_u(t)\ge \rho_u(t)$   for all
$u\in\mmZ^d$ and $t\ge 0$.
\endproclaim

A less obvious property of the construction 
is the following, which forms the basis
of our approach to the hydrodynamic limit. 

\proclaim{Lemma 3.3} Suppose the ballistic deposition
 process
$\sigma$ and a countable family of ballistic
deposition processes $\{\zeta^i:i\in\II\}$
are constructed on a common probability space so that
they all use the same version $\{\TT^u\}$ of the Poisson processes.
 Assume that at time $0$,  almost surely, 
$$\sigma_u(0)=\sup_{i\in\II} \zeta^i_u(0)\quad\text{for all
$u\in\mmZ^d$.}
\tag 3.3
$$
Then almost surely 
$$\sigma_u(t)=\sup_{i\in\II} \zeta^i_u(t)\quad\text{for all
$u\in\mmZ^d$ and $t\ge 0$.}
\tag 3.4
$$
\endproclaim

\demo{Proof} First apply Lemma 3.2 with $\rho=\zeta^i$ 
to get 
$$\sigma_u(t)\ge\sup_{i\in\II} \zeta^i_u(t)\quad\text{for all
$u\in\mmZ^d$ and $t\ge 0$.}
\tag 3.5
$$
Thus we need to show that for all sites $u$ and times $t$,
there is some index $i$ such that 
$\sigma_u(t)= \zeta^i_u(t)$. 

 Pick and fix a realization $\{\TT^u\}$ 
that satisfies assumptions (3.2) and for which the 
conclusion of Lemma 3.1 holds. We first show that, for any
processes $\sigma$ and  $\{\zeta^i:i\in\II\}$ that satisfy
the hypotheses, (3.4) holds for $t\in(0,t_0]$ where $t_0>0$ 
is the number chosen in Lemma 3.1. To do so for a 
fixed site $u^0$, let $\CC\subseteq\mmZ^d$ be the
finite connected component
of $u^0$ in the (random) graph constructed for Lemma 3.1. 
The evolutions of all the processes on the sites of
$\CC$ are determined by the finitely many Poisson 
points in 
 $\cup_{w\in\CC}\TT^w\cap (0,t_0]$. We can now prove 
(3.4) up to time $t_0$ by checking that it 
 holds right after each jump. 

So suppose $r\in (0,t_0]$ is a jump time in $\TT^v$ for some 
site $v\in\CC$ so that 
(3.1) happens for $u=v$. 
Assume
by induction that (3.4) holds for $t< r$, 
 for all $u\in\CC$. Since the variables $\sigma_u$
and $\zeta^i_u$ are $\mmZ$-valued, this means that at each 
time $t< r$ the supremum in (3.4) is actually achieved 
at some $i\in\II$. Depending on how the jump 
(3.1) for $u=v$ is realized, two cases need to be considered.

{\it Case 1.}
First suppose $\sigma_v(r)=\sigma_v(r-)+1$. By induction, there
is a $j\in\II$ such that $\sigma_v(r-)=\zeta^j_v(r-)$. 
Since $\zeta^j_v$ jumps too and by (3.5), 
$$\aligned
\sigma_v(r)&\ge\sup_{i\in\II} \zeta^i_v(r)  \ge\zeta^j_v(r)\\
&= \max\bigl[ \zeta^j_v(r-)+1\,,\,
\max\{ \zeta^j_{v+z}(r-): z\in\NN\}\bigr]\\
&\ge\zeta^j_v(r-)+1 = \sigma_v(r-)+1 \\
&=\sigma_v(r) \,.
\endaligned
$$
Thus $\zeta^j_v$ jumps to the same height as $\sigma_v$. 

{\it Case 2.}
The second possibility is  that for some $w\in\NN$,
$\sigma_v(r)=\sigma_{v+w}(r-)$. Then by the jump rule (3.1)
$$\text{
$\sigma_{v+w}(r-)\ge\sigma_v(r-)+1$ and  
$\sigma_{v+w}(r-)\ge\sigma_{v+z}(r-)$ for all $z\in\NN$.}
$$ 
By the definition of the random graph $v+w\in\CC$,
so by induction there exists a $j\in\II$ 
 such that $\sigma_{v+w}(r-)=\zeta^j_{v+w}(r-)$. 
By (3.5) 
$$\text{
$\sigma_v(r-)+1\ge\zeta^j_v(r-)+1$ and  
$\sigma_{v+z}(r-)\ge\zeta^j_{v+z}(r-)$ for $z\in\NN$.}
$$ 
Together  these equalities and inequalities  
yield 
$$\zeta^j_{v+w}(r-)\ge 
\max\bigl[ \zeta^j_v(r-)+1\,,\,
\max\{ \zeta^j_{v+z}(r-): z\in\NN\}\bigr]\,,
$$
so the jump rule (3.1) applied to $\zeta^j_v$ gives
$$
\zeta^j_v(r)=\zeta^j_{v+w}(r-)=\sigma_{v+w}(r-)=\sigma_v(r).$$
In other words, $\zeta^j_v$ jumped to the same level as
$\sigma_v$, and (3.4) continues to hold right after the jump
time $r$. 

The summarize: we have shown that (3.4) holds for times
$0\le t\le t_0$ for any 
processes $\sigma$, $\{\zeta^i\}$ that satisfy the 
hypothesis (3.3) at $t=0$. Now apply the same step again, to the 
processes $\sigmatil(t)=\sigma(t_0+t)$ and 
$\{\zetatil^i(t)=\zeta^i(t_0+t)\}$.
These processes satisfy the hypothesis at $t=0$ by virtue of (3.4)
at $t=t_0$. This way the validity of (3.4) is extended to 
times $0\le t\le 2t_0$. And so on, to arbitrarily large 
times. 
\qed
\enddemo

For $v\in\mmZ^d$, let $Z^v=(Z^v_u(t):u\in\mmZ^d)$ 
denote the ballistic deposition 
process started from a seed in cell $(v,0)$. In other words,
initially
$$Z^v_u(0)=\cases  0&\text{ if $u=v$}\\
          -\infty&\text{ if $u\ne v$\,.}
\endcases
\tag 3.6
$$
Given an arbitrary 
initial configuration $\sigma(0)=(\sigma_u(0):u\in\mmZ^d)$
[random or deterministic],
define a family of processes, indexed by $\II=\mmZ^d$,
by the initial conditions 
$$\zeta^v_u(0)=\cases  \sigma_v(0)&\text{ if $u=v$}\\
          -\infty&\text{ if $u\ne v$\,.}
\endcases
\tag 3.7
$$
We can write 
$$\zeta^v_u(t)=\sigma_v(0)+Z^v_u(t)
\tag 3.8
$$
with the convention that $\infty+(-\infty)=-\infty$. 
The processes $\sigma$ and $\{\zeta^v:v\in\mmZ^d\}$ 
satisfy (3.3). Lemma 3.3 gives this corollary,
which is basic for our proof of the hydrodynamic limit:

\proclaim{Corollary 3.1} The equality
$$\sigma_u(t)=\sup_{v\in\mmZ^d}\{ \sigma_v(0)
+Z^v_u(t)\}
\tag 3.9
$$ 
holds almost surely, for all $u\in\mmZ^d$ and $t\ge 0$. 
\endproclaim

\hbox{}

\head 4. First-passage site percolation \endhead

In this section we prove some estimates for the 
first-passage percolation problem briefly encountered
in the introduction. First we redefine it in the 
standard way. 

Give each site $u\in\mmZ^d$ an Exp(1)-distributed 
random time  
$t(u)$, independently of the other sites.  Say 
 $\pi=\{w^0,w^1,\ldots,w^m\}\subseteq\mmZ^d$ is a {\it nearest-neighbor
path} from $u$ to $v$ of length $m$ if $m<\infty$,
$w^0=u$, $w^m=v$, and $|w^i-w^{i-1}|=1$ for $i=1,\ldots, m$. 
 We use $|\cdot|$  to denote the $\ell^1$ norm: 
$|w|=|w_1|+\cdots+|w_d|$ for $w=(w_1,w_2,\ldots,w_d)\in\mmZ^d$. 
The passage time  of a path  $\pi=\{w^0,w^1,\ldots,w^m\}$ 
is 
$$T(\pi)=\sum_{i=1}^mt(w^i)\,.
\tag 4.1
$$
Since the path starts from site $w^0$, the value $t(w^0)$
is not included in the sum.  The passage time from
 site $u$ to $v$ is 
$$T(u,v)=\inf_{\pi} T(\pi)
\tag 4.2
$$
where the infimum ranges over nearest-neighbor paths $\pi$
from $u$ to $v$. 
The minimization has the effect that $\pi$ may be
assumed {\it self-avoiding} in the sense that 
there are no repetitions among $\{w^0,w^1,\ldots,w^m\}$. 

The cluster growing from a seed at the origin is defined by 
$$\BB(t)=\{ u\in\mmZ^d: T(0,u)\le t\}\,.
\tag 4.3
$$
It is clear from the description that $\BB(0)=\{0\}$,
and  $\BB(\cdot)$ grows
according to this local rule: Each site adjacent to the 
current cluster joins independently with rate 1. 

To make the connection with ballistic deposition,
consider again the process $Z$ started  from a seed
at the origin:
$$Z_u(0)=\cases 0, &u=0\\
   -\infty, &u\ne 0\,.
\endcases 
\tag 4.4
$$
Let $R(u,h)$ denote the first time $Z$ is at or above 
height $h\in\mmZ_+$ over site $u\in\mmZ^d$:
$$R(u,h)=\inf\{ t>0: Z_u(t)\ge h\}.
\tag 4.5
$$
In particular, $R(u,0)$ is the first time that 
a particle sticks to the cluster above site $u$. 
$R(0,0)=0$ by definition (4.4).
Since the notation may lead to confusion, let us emphasize 
 that the $0$ of $R(u,0)$ is the zero of
$\mmZ_+$, while the $0$ of $T(0,u)$ is the origin of
$\mmZ^d$.

\proclaim{Lemma 4.1}
We have the  following equalities  in 
distribution between the cluster and  passage time  processes: 
$$ \lbrakk \{u\in\mmZ^d: Z_u(t)\ge 0\}\,:\,
t\ge 0\rbrakk   \overset{d}\to=\lbrak \BB(t) : t\ge 0\rbrak 
\tag 4.6
$$
and
$$\{R(u,0):u\in\mmZ^d\}\overset{d}\to=
\{T(0,u):u\in\mmZ^d\}\,.
\tag 4.7
$$
\endproclaim

\demo{Proof}  Abbreviate  $\BBtil(t)=\{u\in\mmZ^d: Z_u(t)\ge 0\}$
for this proof.  
 $\BB(\cdot)$ and $\BBtil(\cdot)$ are both Markov
jump processes on the countable state space 
of finite, connected subsets of $\mmZ^d$ that
contain $0$. They have the same  
 initial state  $\BB(0)=\BBtil(0)=\{0\}$, and
both processes have identical infinitesimal
rates: each site adjacent to the current cluster 
joins independently at rate 1. Hence the two processes are equal in 
distribution, which is (4.6). 

By definition (4.3), 
$T(0,u)$ is the first time when site $u$ joins the cluster
$\BB(t)$.    $R(u,0)$ has the same meaning
for $\BBtil(t)$. So the processes 
$\{T(0,u)\}$ and $\{R(u,0\}$ are obtained by applying
a certain measurable function  to the processes
 $\BB(\cdot)$ and $\BBtil(\cdot)$. 
Thus  (4.7) follows from (4.6). 
\qed
\enddemo

 In later sections 
(4.6) will be used several times to derive 
deviation bounds for $Z$. 

Return to $T(u,v)$ as constructed by (4.3). 
Subadditivity considerations give a limit
$$\lim_{n\to\infty}\frac1n T(0,[nx])=\mu(x)\qquad\text{a.s.}
\tag 4.8
$$
for all $x\in\mmR^d$. The function $\mu$ is
convex, homogeneous [$\mu(rx)=r\mu(x)$ for $r>0$], 
 and Lipschitz continuous. 
 The limiting cluster is defined by
$$B_0=\{x\in\mmR^d: \mu(x)\le 1\}\,,
\tag 4.9
$$
and now the set convergence (2.3) is valid. We shall not
give proofs of these laws of large numbers. 
 The case where the
random times $t(u)$ are on the edges instead of on the 
sites is thoroughly discussed in Kesten's 
lectures \cite{13}. We prove some
large deviation estimates for $T(0,u)$ that we need 
in the sequel. 

\proclaim{Proposition 4.1} For any $x\in\mmR^d$ and $\e>0$ there are
finite  constants
$C_i=C_i(x,\e)>0$ such that 
$$P\bigl( \bigl| T(0,[nx])-n \mu(x) \bigr| \ge n\e\bigr)
\le C_1\exp(-C_2n)
$$
for all $n$. 
\endproclaim

Before the proof we derive two corollaries, assuming 
Proposition 4.1.

\proclaim{Corollary 4.1} For any $t,\e>0$ there are
finite  constants
$C_i=C_i(t,\e)>0$ such that 
$$\sum_{u\in\mmZ^d\,,\,u\notin n(t+\e)B_0}
P\bigl(  T(0,u)\le nt\bigr)
\le C_1\exp(-C_2n)
$$
for all $n$. 
\endproclaim

\demo{Proof of Corollary 4.1} Fix a large $n$
and small $\delta_0,\delta_1,\e_1>0$ so that
$$
\text{$\delta_0+d/n \le \delta_1$ and 
$\delta_1(1+\e_1)\le \e-\e_1$.}
\tag 4.10
$$
 Let $\cAA_m=\bigl[(m+1)(t+\e)B_0
\setminus m(t+\e)B_0\bigr]\cap\mmZ^d$ for $m\ge n$. Write
$$\sum_{u\notin n(t+\e)B_0}
P\bigl(  T(0,u)\le nt\bigr)
= \sum_{m=n}^\infty \sum_{u\in\cAA_m} P\bigl(  T(0,u)\le nt\bigr)\,.
\tag 4.11
$$
Since $|\cAA_m|\le Cm^d$, it suffices to 
bound  $P\bigl(  T(0,u)\le nt\bigr)\le \exp(-Cm)$
for $u\in\cAA_m$. 

Pick and fix points $x^1$, $\ldots$, $x^k$ $\in$
$2(t+\e)B_0\setminus (t+\e)B_0$ such that each 
$y$  $\in$
$2(t+\e)B_0\setminus (t+\e)B_0$ is within $\ell^1$
distance $\delta_0$ of one of the $x^i$'s. 
By the definition (4.9), $\mu(x^i)\ge t+\e$. 
Any $u\in\cAA_m$ satisfies
$$|u-[mx^i]|\le m\delta_0+d
\tag 4.12
$$
for some $1\le i\le k$. For this same $i$,
$$T(0,[mx^i])\le T(0,u)+T(u,[mx^i])\,.$$
$T(u,[mx^i])$ is stochastically dominated 
by a sum $S^1_{[m\delta_1]}$ of $[m\delta_1]$ 
i.i.d.\ Exp(1) random variables, because by (4.12) a direct
lattice path from $u$ to $[mx^i]$ takes at most
 $m\delta_0+d\le m(\delta_0+d/n)
\le m\delta_1$ Exp(1) passage times. Consequently, 
by (4.10), Proposition 4.1, and standard large deviation
bounds for exponential random variables, for $u\in\cAA_m$,
$$\aligned
&P\bigl(  T(0,u)\le nt\bigr)\\
&\le \sum_{i=1}^k P\bigl(  T(0,[mx^i])\le m\mu(x^i)-m\e_1\bigr)
+ k P\bigl( S^1_{[m\delta_1]}\ge m\delta_1(1+\e_1)\bigr)\\
&\le C_1k\exp(-C_2m)
\endaligned
$$
for some finite constants $C_1, C_2>0$ [not the same 
as those in Proposition 4.1]. By substituting  this bound 
in (4.11) we have proved Corollary 4.1
for large enough $n$. By increasing $C_1$ in the 
statement it follows for all $n$. 
\qed
\enddemo

\proclaim{Corollary 4.2} For any $t,\e>0$ there are
finite  constants
$C_i=C_i(t,\e)>0$ such that 
$$ P\bigl(  n(t-\e)B_0\nsubseteq \BB(nt)\bigr)
\le C_1\exp(-C_2n)
$$
for all $n$. 
\endproclaim

\demo{Proof} The event $n(t-\e)B_0\nsubseteq \BB(nt)$
implies that $T(0,u)>nt$ for some $u\in n(t-\e)B_0$.
The idea for proving this is the same as for the 
previous Corollary: Pick a fine grid of point
 $x^1$, $\ldots$, $x^k$ $\in$
$(t-\e)B_0$. For each  $u\in \bigl[n(t-\e)B_0\bigr]\cap\mmZ^d$ pick
the closest $[nx^i]$, use the inequality
$$T(0,u)\le T(0,[nx^i])+T([nx^i],u)\,,$$
dominate $T([nx^i],u)$ by a sum of i.i.d.\ Exp(1)'s,
 note that $\mu(x^i)\le t-\e$, and apply Proposition 4.1.
We leave the details to the reader. 
\qed
\enddemo

As the last item of this section we prove Proposition 4.1.
For this purpose we combine a result of Talagrand \cite{30}, 
stated as the next lemma, with some ideas from
Kesten \cite{14}. Suppose $\{X_i:1\le i\le N\}$ are 
independent random variables such that $0\le X_i\le 1$.
Let $\FF$ be a family of $N$-tuples $\a=(\a_i:1\le i\le N)$
of real numbers.
Set $\sigma=\sup_{\a\in\FF}\|\a\|_2$, where 
$\|\a\|_2=\biggl(\sum_{i=1}^N \a_i^2\biggr)^{1/2}$. 
Define the random variable 
$$Z=\sup_{\a\in\FF}\sum_{i=1}^N \a_iX_i\,.
\tag 4.13
$$

\proclaim{Lemma 4.2} {\rm \cite{30}} For any numbers $a<b$, 
$$P\bigl( Z\ge b\bigr)\cdot P\bigl( Z\le a\bigr)
\le \exp\bigl[-(b-a)^2/4\sigma^2\bigr]\,.
\tag 4.14
$$
\endproclaim

Some comments: We use here the original notation of \cite{30}
 even though the use of $Z$, $\sigma$ and $\FF$ 
 conflicts with our notation elsewhere.
 But this notation is used only for the proof of Proposition 4.1  
so it will not cause confusion. 
Except for Lemma 4.2 our proof of Proposition 4.1 is
 self-contained, and the  proof of Lemma 4.2 can be 
quickly read from  Talagrand's paper \cite{30}: 
 the reader only needs to 
go through section 4.1 which is self-contained, and  the
proof of Theorem 8.1.1 which relies on section 4.1. 
Our Lemma 4.2 is contained in the proof of 
Theorem 8.1.1 in \cite{30}. 

Fix $x\in\mmR^d$. 
To produce the random variable $Z$ of (4.13), we truncate
the random times of the sites and restrict the set 
of paths $\pi$. Let $K, L>0$ be constants independent 
of $n$, to be chosen below. Let 
$\Pi_{L}$ denote the collection of self-avoiding
nearest-neighbor paths $\pi$ from $0$ to $[nx]$ with 
length at most $Ln$: 
$\pi=\{0=w^0, w^1,\ldots, w^{m-1}, w^m=[nx]\}$ such that $m\le Ln$. 
Define truncated times by $\that(u)=t(u)\wedge K$ for $u\in\mmZ^d$. 
Let $\That$ denote the passage time with truncated
variables and restricted paths:
$$\That(0,[nx])=\min_{\pi\in\Pi_{L}}\That(\pi)\,,
$$
where $\That(\pi)$ is defined by (4.1) with $t$ replaced by $\that$.

 Let 
$\{0=u^0, u^1, u^2,\ldots , u^N\}$ 
be a numbering of the sites  that are within $\ell^1$
distance $Ln$ of the origin.
Let  $\FF$ be the subset of $\{-1,0\}$-valued 
 $N$-tuples 
$\a=(\a_i)$ that are negatives of indicator
functions of paths $\pi\in\Pi_L$: $\a\in\FF$ iff for
some $\pi\in\Pi_L$,   $\a_i=-\ind\{u^i\in\pi\}$ for all $1\le i\le N$. 
   Define i.i.d.\ random
variables with values in $[0,1]$ by $X_i=K^{-1}\that(u^i)$. The
 random variable $Z$ defined in (4.13) is then 
$$Z=\sup_{\a\in\FF}\sum_{i=1}^N\a_i\frac{\that(u^i)}K=
-\frac1{K}\inf_{\pi\in\Pi_L}\sum_{i: u^i\in\pi}\that(u^i)=
-\frac1{K}\That(0,[nx])\,.
\tag 4.15
$$

\proclaim{Lemma 4.3} Given $\e>0$, we can choose the  constants
$L=L(x)$ and $K=K(\e,x)$ so that 
$$P\bigl(\, \bigl| \That(0,[nx])- T(0,[nx]) \bigr|
\ge n\e \bigr) \le C_1\exp(-C_2n)
\tag 4.16
$$
for all $n$,  for some finite 
constants $C_1, C_2>0$. 
\endproclaim

\demo{Proof} $T(0,[nx])$ has a.s.\ a unique optimal 
path $\pi_n$ because the site times $t(u)$ have a 
continuous distribution. $\That(0,[nx])$ does not necessarily have a
unique optimal path, so order the paths in $\Pi_L$ in some way
and let  $\pihat_n$ denote the first path in this ordering
that is  optimal for $\That(0,[nx])$. 
 Divide the estimation in (4.16) into two parts:  
$$\aligned
&P\bigl(\, \bigl| \That(0,[nx])-T(0,[nx])\bigr|
\ge n\e \bigr)\\
&\le  P\bigl( \pi_n\notin \Pi_{L}\bigr)
 + P\bigl(\, \bigl| \That(0,[nx])- T(0,[nx])\bigr|
\ge n\e  \,,\,\pi_n\in \Pi_{L} \bigr)\,.
\endaligned
$$
The first term of the right-hand side is bounded
as follows: 
$$\aligned
&P\bigl( \pi_n\notin \Pi_{L}\bigr) \\
&\le P\bigl(\,  T(0,[nx]) \ge nL_0  \bigr)
 + P\biggl(\text{ a self-avoiding path 
$\pi$ starts from $0$}\\
&\qquad\qquad \text{ and has length at least $Ln$ but
$T(\pi)<nL_0$ }\biggr)\\
&\le P\bigl(\, S^1_{|[nx]|} \ge nL_0  \bigr)
 + \sum_{m=Ln}^\infty (2d-1)^m P\biggl( S^1_m < m \cdot
\frac{L_0}{L}\,\biggr)\,.
\endaligned
$$
Above we used $T(0,[nx])\le S^1_{|[nx]|}$, 
 bounded the number of self-avoiding
paths from $0$ with length $m$ by $(2d-1)^m$ [at each step at most
$2d-1$ sites to choose from], and noted that 
for a path $\pi$ of length $m$, $T(\pi)$ is distributed
like $S^1_m$,  a sum
of $m$ i.i.d.\ Exp(1) variables. By choosing
$L_0 = x+1$ and  $L=L(x)=2x+2$
the quantity above is bounded by 
$$ \exp(-C_1n)+\sum_{m= Ln}^\infty (2d-1)^m\exp(-C_2m)\le C_1\exp(-C_2n)$$ 
for 
constants $C_1,C_2>0$. Note that the constants change from 
one inequality to the next.

To bound the probability $P\bigl(\, \bigl| \That(0,[nx])- T(0,[nx])\bigr|
\ge n\e  \,,\,\pi_n\in \Pi_{L} \bigr)$
note first that if $\pi_n\in \Pi_{L}$, then
necessarily 
$$\aligned
0&\le T(0,[nx])- \That(0,[nx])
\le \sum_{w\in\pihat_n} \bigl[ t(w)-\that(w)\bigr]\\
&= \sum_{w\in\pihat_n} [ t(w)-K]\cdot\ind\{t(w)>K\}\,.
\endaligned
$$
For paths $\pi\in\Pi_{L}$ define the events
$$A_\pi=\lbrakk \,\sum_{w\in\pi} [ t(w)-K]\cdot\ind\{t(w)>K\}
\ge n\e  \rbrakk
\quad\text{and }\quad
D_\pi=\{\pihat_n=\pi\}\,.$$
Let $\GG_\pi$ denote the 
$\sigma$-algebra generated by $\{t(w): w\notin \pi\}$. 
When $\pi$ and  $\{t(w): w\notin \pi\}$ are fixed, 
$A_\pi$ is an increasing event and $D_\pi$ is a decreasing
event of the variables  $\{t(w): w\in \pi\}$. 
We estimate as follows:
$$\aligned
&P\bigl(\, \bigl| \That(0,[nx])- T(0,[nx])\bigr|
\ge n\e  \,,\,\pi_n\in \Pi_{L} \bigr)\\
&\le 
P\biggl(\sum_{w\in\pihat_n} [ t(w)-K]\cdot\ind\{t(w)>K\}\ge n\e  
\biggr)\\
&=\sum_{\pi\in\Pi_{L}}P\bigl( A_\pi\cap D_\pi\bigr)\\
&=\sum_{\pi\in\Pi_{L}}E\bigl[ P\bigl( A_\pi\cap D_\pi \bigm|\GG_\pi\bigr)
\bigr]\\
&\le \sum_{\pi\in\Pi_{L}}E\bigl[ 
P( A_\pi|\GG_\pi) P( D_\pi|\GG_\pi)\bigr]\\
&\qquad\qquad\qquad\qquad[\text{ by the FKG inequality }]\\
&= \sum_{\pi\in\Pi_{L}} P( A_\pi)\,E\bigl[ 
 P( D_\pi|\GG_\pi)\bigr]\\
&\qquad\qquad\qquad\qquad[\text{ $A_\pi$ is independent of $\GG_\pi$ }]\\
&\le P\biggl(\,\sum_{i=1}^{Ln} [ t(u^i)-K]\cdot\ind\{t(u^i)>K\}\ge n\e 
\biggr)\cdot
\sum_{\pi\in\Pi_{L}} P( D_\pi)\\
&\qquad\qquad\qquad\qquad[\text{ each $\pi$ has at most $Ln$ 
passage times }]\\
&\le P\biggl(\,\sum_{i=1}^{Ln} [ t(u^i)-K]\cdot\ind\{t(u^i)>K\}\ge n\e 
\biggr)\,.
\endaligned
$$
This last probability is $\le \exp(-Cn)$ if we
choose $K=K(\e,L)$ large enough so that 
$E[(t(u^i)-K)\cdot\ind\{t(u^i)>K\}]<\e/(2L)$. 
\qed
\enddemo

We are ready to prove Proposition 4.1. For convenience,
replace $\e$ by $8\e$ in the statement of the 
Proposition. Fix $\e>0$, and choose
 $K,L$  so that (4.16) holds.  
$$\aligned
&P\bigl(\, \bigl| T(0,[nx])-n\mu(x)\bigr|
\ge 8n\e \bigr)\\
&\le 
P\bigl(\, \bigl| T(0,[nx])-\That(0,[nx])\bigr|
\ge 4n\e \bigr)\\
&\qquad +
P\bigl(\, \bigl| \That(0,[nx])-n\mu(x)\bigr|
\ge 4n\e \bigr)\,,
\endaligned
$$
so by Lemma 4.3 it suffices to bound the last probability.
By the limit in (4.8) and by Lemma 4.3, 
$$P\bigl(\, \bigl| \That(0,[nx])-n\mu(x)\bigr|
\ge 2n\e \bigr)< 1/2
\tag 4.17
$$
 for all large enough $n$. 
Recall equality (4.15), use (4.17)  and estimate as follows: 
$$\aligned
&P\bigl(\, \bigl| \That(0,[nx])-n\mu(x)\bigr|
\ge 4n\e \bigr)\\
&=P\bigl(\, \That(0,[nx])\ge n\mu(x) + 4n\e \bigr)\\
&\qquad +\ P\bigl(\, \That(0,[nx])\le n\mu(x) - 4n\e \bigr)\\
&=P\bigl(\, Z\le -n\mu(x)/K - 4n\e/K \bigr)\\
&\qquad +\ P\bigl(\, Z\ge -n\mu(x)/K + 4n\e/K \bigr)\\
&\le 2\cdot P\bigl(\, Z\le -n\mu(x)/K - 4n\e/K \bigr)
\cdot P\bigl(\, Z\ge -n\mu(x)/K - 2n\e/K \bigr)\\
&\qquad +\ 2\cdot P\bigl(\, Z\ge -n\mu(x)/K + 4n\e/K \bigr)
\cdot P\bigl(\, Z\le -n\mu(x)/K + 2n\e/K \bigr)\\
&\le 4\exp\bigl[-\bigl(2n\e/K\bigr)^2/4Ln\bigr]\\
&=4\exp\bigl[-\e^2n/(K^2L)\bigr]\,.
\endaligned
$$
We used Lemma 4.2 in the second last step. Each $\a\in\FF$ has
 at most $Ln$ entries equal to $-1$
and the rest are zeroes. Consequently $\sigma\le\sqrt{Ln\,}$. 
We have proved Proposition 4.1 for large enough $n$, and it follows
for all $n$ by increasing $C_1$ sufficiently.

\hbox{}

\head 5. The ballistic deposition shape from a seed \endhead

In this section we prove Theorem 1, the almost sure limit for the 
ballistic deposition process started from a single
seed. Initially 
$$Z_u(0)=\cases  0&\text{ if $u=0$}\\
          -\infty&\text{ if $u\ne 0$\,.}
\endcases
\tag 5.1
$$
Caution to the reader: $Z$ was defined by (4.15) only
for the proof of Proposition 4.1. Elsewhere in
the paper, $Z$ denotes the ballistic deposition
process started from a seed. 
The process $Z$ is constructed by the argument of Section 3
on the probability space $(\Omega, \FF,P)$ of the Poisson
jump time processes $\TT=\{\TT^u\}$. 
 Recall that 
$B_0\subseteq\mmR^d$ denotes
the closed, convex limiting set for first-passage site percolation 
on $\mmZ^d$ with Exp($1$) waiting times. The goal is 
to prove that for a bounded, positive, concave function
$g$ defined on the open set $\inte B_0$, 
$$\lim_{n\to\infty} \frac1n Z_{[nx]}(nt)=tg(x/t)\qquad\text{a.s.}
\tag 5.2
$$
for all $x\in\mmR^d$ and $t>0$ such that $x/t\in\inte B_0$.

The main tool in the proof is the Kesten-Hammersley lemma
from subadditive ergodic theory.
Since $Z$ is unbounded both above and below, we
work instead with the passage times $R(u,h)$ 
defined by (4.5).
By definition $R(u,h)\ge 0$. Since $Z$ can always
reach cell $(u,h)$ by advancing along
each coordinate axis in turn, with Exp(1) waiting times,
$$R(u,h) \sle S^1_{|u|+h}\,.
\tag 5.3
$$
Here $\sle$ denotes stochastic dominance
and $S^1_{n}$ is a sum of $n$ i.i.d.\ Exp($1$) random
variables. Recall that $|u|$ is the $\ell^1$ norm on $\mmZ^d$.
In particular, $R(u,h)$ has finite moments of all order.

The main property to check is a subadditivity:

\proclaim{Lemma 5.1}
For $u,v\in\mmZ^d$ and $h,k\in\mmZ_+$, 
there exists  a random variable $\Rtil(v,k)$ such that
$$R(u+v,h+k)\le R(u,h)+\Rtil(v,k)\,,
\tag 5.4
$$
and $\Rtil(v,k)$ is 
independent of $R(u,h)$, and equal in distribution to
$R(v,k)$. 
\endproclaim

\demo{Proof}
To define  $\Rtil(v,k)$, we start a new $Z$-process
at time $R(u,h)$ from a seed
in cell $(u,h)$. This new process 
$\Ztil$ is dominated by the original $Z$-process,
hence the inequality (5.4). 

For the reader not familiar with these types of arguments,
here are the details: Think of  $\TT=\{\TT^w:w\in\mmZ^d\}$
as an infinite-dimensional vector of Poisson point processes  $\TT^w$,
 indexed by time $t\in(0,\infty)$. 
  Let $\FF_t$ be the $\sigma$-algebra generated
by the restriction of $\TT$ to the time-interval $(0,t]$.
Then $R(u,h)$ is a stopping
time for the filtration $(\FF_t)$.
Restart the Poisson processes at time $R(u,h)$ and 
translate the index by $u$ to get new point processes
$\TTtil=\{\TTtil^w:w\in\mmZ^d\}$ where 
$\TTtil^w=\bigl[\TT^{u+w}-R(u,h)\bigr]\cap(0,\infty)$. 
(The subtraction means that
epochs of $\TT^{u+w}$ are translated back $R(u,h)$ time units.)

\proclaim{Sublemma} $\TTtil$ is a collection of i.i.d.\ rate 1
Poisson processes on $(0,\infty)$, independent of $R(u,h)$. 
\endproclaim

\demo{Proof of sublemma} We ignore the spatial translation $w\mapsto u+w$. 
$R=R(u,h)$ is the stopping time, and $\FF_R$ is the 
$\sigma$-algebra up to time $R$, consisting of events 
$A$ such that $A\cap\{R\le t\}\in\FF_t$ for all $t$. 
 Topologize the 
space of families  $\TT=\{\TT^w:w\in\mmZ^d\}$
of point processes on $(0,\infty)$ by
the product of the vague topology on point processes.
Let $g$ be a bounded
continuous function on that space. 
 The claim is that for any $A\in\FF_R$ and any such
 $g$, 
$$E[\ind_A\, g(\TTtil)]=P(A)E[g(\TT)]\,.
\tag 5.5
$$ 
Since $R$ itself is $\FF_R$-measurable, the lemma follows. 

 Take a sequence of discrete stopping
times $R_n$ that decrease down to $R$, almost surely. 
For example,  $R_n=\sum_k k2^{-n}\ind\{(k-1)2^{-n}<R
\le k2^{-n}\}.$ 
Let $\TTtil_n$ denote $\TT$ 
translated back $R_n$ time units and restricted to $(0,\infty)$.
 $A\in\FF_R$ implies  $A\in\FF_{R_n}$,
and we get: 
$$
\aligned
E[\ind_A\, g(\TTtil_n)]&=\sum_{i=1}^\infty 
E[\ind_A\, \ind\{R_n=t_i\} g(\TTtil_n)]\\
&=\sum_{i=1}^\infty 
E\bigl[\ind_A \ind\{R_n=t_i\} g\bigl((\TT-t_i)\cap(0,\infty)\bigr)\bigr]\\
&=\sum_{i=1}^\infty 
E[\ind_A \ind\{R_n=t_i\}] E[ g(\TT)]\\
&=P(A)E[g(\TT)]\,.
\endaligned
$$ 
The $\{t_i\}$ are the possible values of $R_n$. The point
of the above calculation 
is of course that, when $t_i$ 
is a deterministic number, $A\cap\{R_n=t_i\}\in\FF_{t_i}$,
while $\TT-t_i$ restricted to $(0,\infty)$
 is an i.i.d.\ collection of Poisson 
point processes independent of $\FF_{t_i}$. Now let $n\to\infty$, and 
observe that $\TTtil_n\to\TTtil$ a.s.\ so that 
$E[\ind_A\, g(\TTtil_n)]$ converges to $E[\ind_A\, g(\TTtil)]$,
and we obtain (5.5) in the limit. 
\qed
\enddemo

Return to the proof of Lemma 5.1. 
Define the process $\Ztil$ as a function of $\TTtil$,
exactly as $Z$ is a function of $\TT$, with an initial
seed at the origin: $\Ztil_w(0)=-\infty\cdot\ind\{w\ne 0\}$
as in (5.1) for $Z$. Let $\Rtil(v,k)$ be the time when $\Ztil$ gets
at or above cell $(v,k)$. Then $\Rtil(v,k)$ is independent
of $R(u,h)$, and has exactly the same distribution
as $R(v,k)$ defined by (4.5).

  Consider the processes
$\sigma$ and $\rho$ defined by 
$$\text{$\sigma_w(t)=Z_w(R(u,h)+t)$ and 
$\rho_w(t)=h+\Ztil_{w-u}(t)$.}
$$
 Then $\rho_w(0)=-\infty$
except at $w=u$, where 
$$\rho_u(0)=h+\Ztil_{0}(0)=h=
Z_u(R(u,h))=\sigma_u(0).
$$
 The jump times
of $\sigma_w$ are given by $\TT^w-R(u,h)$, and
those of $\rho_w$ by  $\TTtil^{w-u}=\TT^w-R(u,h)$. 
So processes $\sigma$ and $\rho$ use the same 
Poisson jump times and initially 
 $\sigma(0)\ge \rho(0)$. By 
Lemma 3.2 $\sigma(t)\ge \rho(t)$ for all
$t$. Take $t=\Rtil(v,k)$. Then 
$$\aligned
Z_{u+v}(R(u,h)+\Rtil(v,k))&=
\sigma_{u+v}(\Rtil(v,k))\\
&\ge \rho_{u+v}(\Rtil(v,k))
=h+\Ztil_{v}(\Rtil(v,k))\\
&=h+k\,,
\endaligned$$
which says that $Z_{u+v}$ has reached height $h+k$ by
time $R(u,h)+\Rtil(v,k)$. This   implies (5.4) and 
completes the proof of Lemma 5.1. 
\qed
\enddemo 

  Inequality (5.4), the existence of moments, and 
the Kesten-Hammersley lemma as given on p.\ 20 of 
\cite{28} imply that for all $u\in\mmZ^d$ and $h\in\mmZ_+$
there exists a number $\gamma(u,h)$ such that, for
any positive integer $m$, 
$$\lim_{j\to\infty} \frac1{2^jm} R(2^jmu, 2^jmh) =\gamma(u,h)
\qquad \text{a.s.}
\tag 5.6
$$

The first task is to extend (5.6) to a genuine limit, and for
that we use the following
 continuity property of the passage times. It is an
immediate consequence of (5.3), (5.4) and the monotonicity
of $R(u,h)$ in the $h$-variable. 

\proclaim{Lemma 5.2} For any $u,v\in\mmZ^d$ and $h,k\in\mmZ_+$,
$$R(u,h)-R(v,k)\sle S^1_{|u-v|+(h-k)^+}\,.
\tag 5.7
$$
\endproclaim

Now regard $(u,h)$ fixed, and also
fix $\e>0$ and an integer $m$ large enough relative 
to $\e$ and $(u,h)$ [how large $m$ needs to be is seen shortly]. 
For large enough $n$, pick first $j=j(n)$
so that $2^jm\le n< 2^{j+1}m$. Then it is possible 
to find a $k=k(n)\in\{0,1,\ldots,m-1\}$ such that
$$2^j(m+k)\le n< 2^{j}(m+k+1).
\tag 5.8
$$
 We write 
$$\aligned
&\frac{2^j(m+k)}n\cdot\frac1{2^j(m+k)}
R\bigl(2^j(m+k)u,2^j(m+k)h\bigr)\\
&\le \frac1n R(nu,nh) +\frac1n U(n) \,,
\endaligned
\tag 5.9
$$
where the error
$$U(n)\equiv \lbrak R\bigl(2^j(m+k)u,2^j(m+k)h\bigr)-R(nu,nh)
\rbrak \vee 0 
\sle S^1_{2^jL}$$
 for a constant $L=L(u,h)$, by (5.7). 
The Cram\'er rate function for Exp(1) is 
$\kappa(x)=x-1-\log x$, so we get the estimate
$$\aligned
P\bigl( U(n) \ge n\e\bigr) &\le P\biggl(S^1_{2^jL} 
 \ge 2^jL\frac{n\e}{2^jL}\biggr)  \\
&\le P\biggl(S^1_{2^jL} 
 \ge 2^jL m\e/L\biggr)  \\
&\le \exp\bigl[-2^jL \,\kappa\bigl( m\e/L\bigr)\bigr]\\
&\le \exp[-Cn]\,,
\endaligned
\tag 5.10
$$
where $C>0$ is a constant. In the above calculation
  we used $2^jm\le n< 2^{j+1}m$ and $\kappa(x)\ge x/2$ 
for large enough $x$, and took $m$ large enough.

Let $n\to\infty$ in (5.9), so that $j\to\infty$ also.
Even though $k=k(n)$ varies with $n$, it has only finitely
many possible values so the limit (5.6) happens 
 on the left-hand side of
 (5.9). Note  that
$$\frac{2^j(m+k)}n\ge \frac{m+k}{m+k+1}\ge \frac{m}{m+1}\,.$$
 The error
$n^{-1}U(n)$ vanishes a.s.\ by   the estimate  (5.10) and 
 Borel-Cantelli. We get 
$$\liminf_{n\to\infty} \frac1n R(nu,nh)\ge
\frac{m}{m+1} \gamma(u,h)\qquad
\text{a.s.}$$
A similar argument works for the limsup. 
Let $m\to\infty$, and we  
 have this 
intermediate statement:  for all $u\in\mmZ^d$ and $h\in\mmZ_+$
there exists a number $\gamma(u,h)$ such that
$$\lim_{n\to\infty} \frac1n R(nu,nh)= \gamma(u,h)\qquad
\text{a.s.}
\tag 5.11
$$

From (5.7) we get a Lipschitz property for $\gamma$:
$$|\gamma(u,h)-\gamma(v,k)|\le |u-v|+|h-k|\,,
\tag 5.12
$$
while (5.11) gives  homogeneity: for nonnegative integers $m$
$$\gamma(mu,mh)=m\gamma(u,h)\,.
\tag 5.13
$$
 (5.13) permits us to
define $\gamma$ unambiguously for  $(x,b)\in\mmQ^d\times\mmQ_+$ by 
$$ \gamma(x,b)=\frac1m \gamma(mx,mb)
\tag 5.14
$$ where 
$m$ is any positive integer such that 
$(mx,mb)\in\mmZ^d\times\mmZ_+$. 
By an estimation similar to that used in (5.9), the limit
in (5.11) can be extended to all $(u,h)\in\mmQ^d\times\mmQ_+$.

The final step is to extend the limit in (5.11) 
so that it holds outside a single exceptional
$P$-null set for all  $(u,h)=(x_0,b_0)\in\mmR^d\times\mmR_+$.
The Lipschitz property (5.12) continues to hold on 
$\mmQ^d\times\mmQ_+$ for the extension of $\gamma$
defined by (5.13), so we can extend $\gamma$ uniquely
to a Lipschitz function on $\mmR^d\times\mmR_+$. 
For rational $(x,b)$, integers $n$, 
and rational $\delta>0$ define the event $A_{n,\delta}(x,b)$
by
$$\aligned
A_{n,\delta}(x,b)=\lbrakk &\text{there exists 
$(v,k)\in\mmZ^d\times\mmZ_+$ such that}\\
 &\text{$\bigl|[nx]-v\bigr|+\bigl|[nb]-k\bigr|\le n\delta$ but  }\\
 &\text{$\bigl|R\bigl([nx],[nb]\bigr)-R(v,k)\bigr|\ge 2n\delta$}
\rbrakk\,.
\endaligned
\tag 5.15
$$
By (5.7) and by standard large deviation bounds for exponential
random variables, 
$$P\bigl(A_{n,\delta}(x,b)\bigr)\le 
C_1n^{d+1}\exp[-C_2n] 
$$
for finite constants $C_i=C_i(x,b,\delta)>0$. 
This bound is suitable for Borel-Cantelli. Thus at this
stage the following holds with 
probability 1: for each $(x,b)\in\mmQ^d\times\mmQ_+$
and rational $\delta>0$, (5.11) holds 
with $(u,h)=(x,b)$, and for large enough $n$, 
$$
\bigl|R\bigl([nx],[nb]\bigr)-R(v,k)\bigr|\le 2n\delta
\tag 5.16
$$
for all cells $(v,k)$ that satisfy
$$\bigl|[nx]-v\bigr|+\bigl|[nb]-k\bigr|\le n\delta.
\tag 5.17
$$ 

Now let $(x_0,b_0)\in\mmR^d\times\mmR_+$. Choose rational
$(x,b)$ such that $|x-x_0|+|b-b_0|\le \delta/2$. Take
 $(v,k)=([nx_0],[nb_0])$.   Then 
(5.17) holds for large enough $n$,
and by letting $n\to\infty$ in (5.16) we get 
$$\gamma(x,b)-2\delta \le \liminf_{n\to\infty}\frac1n 
 R([nx_0],[nb_0])\le \limsup_{n\to\infty}\frac1n 
 R([nx_0],[nb_0])\le \gamma(x,b)+2\delta.
$$
Take rational $(x,b)$ that converge to $(x_0,b_0)$, and 
use the Lipschitz continuity of $\gamma$. This 
proves the limit in this proposition:

\proclaim{Proposition 5.1}  There is a homogeneous,
subadditive, convex
Lipschitz function $\gamma$ on  $\mmR^d\times\mmR_+$
such that, outside an event of probability zero,
$$\lim_{n\to\infty} \frac1n R([nx],[nb])= \gamma(x,b)
\tag 5.18
$$
for all $(x,b)\in\mmR^d\times\mmR_+$.
\endproclaim

To prove Proposition 5.1, it remains 
to argue the properties of $\gamma$: The 
subadditivity (5.4) implies the corresponding subadditivity
for $\gamma$. This subadditivity is preserved by the
extensions of $\gamma$ to rational and real points 
of $\mmR^d\times\mmR_+$. Same is true of the homogeneity
(5.13). Homogeneity and subadditivity together imply convexity. 

At this point we have not ruled out the possibility
that $\gamma\equiv 0$. (5.3) gives the upper bound
$$\gamma(x,b)\le |x|+b\qquad\text{for all $(x,b)\in\mmR^d\times\mmR_+$. }
\tag 5.19
$$
To get a positive lower bound for $\gamma$, we construct
another growth process in $\mmZ^d\times\mmZ_+$ whose 
height dominates $Z$, and that has a simple structure so that its
spread is easier to bound. 

Instead of just focusing on the height $Z_u(t)$ of the 
growing cluster, let us denote by $\ZZ(t)$ the actual set of
occupied cells at time $t$. $\ZZ(t)$ is a subset of 
 $\mmZ^d\times\mmZ_+$, and initially $\ZZ(0)=\{(0,0)\}$. 
The rule of evolution for $\ZZ(\cdot)$ is this:
At epochs of $\TT^u$, the top growth cell above $u$
is annexed to $\ZZ$. By definition, the top growth cell
above site $u$ is $(u,k)\in\mmZ^d\times\mmZ_+$
with maximal $k$ subject to the condition that 
$$
\text{$(u,k)\notin\ZZ$, but
either $(u,k-1)\in\ZZ$ or $(u+z,k)\in\ZZ$ for some
$z\in\NN$.}
\tag 5.20
$$
If no finite $k$ satisfies this condition, nothing is
annexed to $\ZZ$.  

Define a different growing cluster $\WW$ by stipulating that {\it all}
growth cells (and not just the top one) join independently
with rate 1. A cell $(u,k)$ is a growth cell for $\WW$
if (5.20) holds with $\WW$ instead of $\ZZ$. So note in particular
that the cluster $\WW$ grows ``sideways'' and ``up''
in $\mmZ^d\times\mmZ_+$, but not down.
Set  initially $\WW(0)=\{(0,0)\}$. 

To construct $\WW$ we employ a collection $\{\TT^u_q: u\in\mmZ^d, 
q\in\mmZ_+\}$ of i.i.d.\ Poisson point processes. The nonnegative 
integer $q$ labels the growth cells from top down, so that
the top growth cell above site $u$
 is assigned label $q=0$, the next
highest growth cell gets label $q=1$, and so on. More 
formally, we can define 
$$g_{u,0}=\max\{ k: \text{$(u,k)$ is a growth cell for $\WW$}\}
$$
(``$g$'' for growth) and inductively for $q\ge 1$
$$g_{u,q}=\max\{0\le  k<g_{u,q-1} : 
\text{$(u,k)$ is a growth cell for $\WW$}\}\,.
$$
By convention, the maximum of an empty set is $-\infty$, so 
if the current cluster $\WW$ has exactly $m$ growth cells above
site $u$, then 
$g_{u,q}=-\infty$ for $q\ge m$. The $g_{u,q}$'s are of
course functions of time too. 
The precise rule of evolution for $\WW$ is this: 
$$
\aligned
&\text{if $r$ is an epoch of $\TT^u_q$ and $g_{u,q}>-\infty$, then
$\WW(r)=\WW(r-)\cup\{(u,g_{u,q})\}$.}\\
&\text{Update the $g_{v,q}$'s for $v=u$ and $v\in u+\NN$.}
\endaligned
\tag 5.21
$$

The argument given for the construction of ballistic deposition
in Section 3 does not work for $\WW$ because each site $u$ now
has infinitely many Poisson processes $\{\TT^u_q:q\in\mmZ_+\}$ attached 
to it. However, we can easily see that, given any finite
time $t_1$, at most finitely many $\TT^u_q$-processes are
involved in constructing the dynamics up to time $t_1$: 
Starting with the seed at $(0,0)$, let $\tau_1$, $\tau_2$,
$\tau_3$, $\ldots$ be the successive waiting times for 
adding the second, third, fourth,...\ cell to the 
existing cluster. Since each new particle adds at most 
$2d+1$ growth cells, $\tau_n$ is stochastically larger 
than an Exp($2dn+n$) random variable. 
Consequently $\sum_n\tau_n=\infty$ a.s.\ and 
only finitely many steps are needed (and only finitely
many Poisson processes $\TT^u_q$ inspected) to construct
$\WW(t)$ for $0\le t\le t_1$. 

Now couple $\WW(t)$ and $\ZZ(t)$ by letting $\{\TT^u_0\}$ be
the Poisson processes that govern the evolution of $\ZZ$. 
In other words, both $\ZZ$ and $\WW$ annex their top 
growth cell above $u$ simultaneously at epochs of $\TT^u_0$. 
At epochs of $\TT^u_q$ for $q\ge 1$ the $\ZZ$-process
 does nothing, 
while $\WW$ may add other growth cells as stipulated in
rule (5.21). 
Since the top growth cells are not necessarily
 the same for $\ZZ$ and
$\WW$ it does not follow that $\WW(t)$ contains $\ZZ(t)$,
but it does follow that the height of $\WW$ always
dominates the height of $\ZZ$. If we let 
$$W_u(t)=\max \{k\ge 0: (u,k)\in\WW(t)\} $$
when some $(u,k)$ lies in $\WW(t)$ and $W_u(t)=-\infty$
otherwise, we get this inequality:
$$Z_u(t)\le W_u(t)\qquad \text{for all $u\in\mmZ^d$ and $t\ge 0$, a.s.}
\tag 5.22
$$
Proof of (5.22) is by induction on jumps [which  
are only finitely many in any finite time interval, almost
surely, as argued above]. As long as (5.22) holds, no 
top growth cell of $\ZZ$ can be  above the corresponding
 top growth
cell of $\WW$, and consequently the next jump 
preserves (5.22).  

To make use of (5.22) we redefine $\WW$ as a first-passage
problem. Give the cells  i.i.d.\ Exp($1$) random waiting times $\{t(u,h): 
(u,h)\in\mmZ^d\times\mmZ_+\}$. Consider 
 self-avoiding nearest-neighbor
lattice paths $\pi=\{(v^0,k^0), (v^1,k^1), \ldots, (v^m,k^m)\}$
whose  admissible
steps are of these types: for each $i=1,\ldots,m$,  
$$\text{ $(v^i,k^i)- (v^{i-1},k^{i-1}) = (0,1)$ or 
$(\pm \ehat_p,0)$ for some $p=1,\ldots,d$,}
\tag 5.23
$$
 where $\ehat_1$,
$\ldots$, $\ehat_d$ are the $d$ standard basis vectors
in $\mmR^d$. In other words, inside a layer 
$\mmZ^d\times\{k\}$ an admissible path $\pi$ takes arbitrary 
nearest-neighbor steps subject to self-avoidance, 
and across the layers $\pi$ moves only up, not down.
The admissible steps are chosen to match (5.20). 
 The passage time 
of such a path $\pi$ is $M(\pi)=\sum_{i=1}^mt(v^i,k^i)$. The 
passage time of cell $(u,h)$ is 
$$M(u,h)=\inf_{\pi} M(\pi) 
$$
where the infimum ranges over paths  $\pi$ of the above
type from  $(0,0)=(v^0,k^0)$ to  $(v^m,k^m)=(u,h)$. Let 
$$\WWtil(t)=\{(u,h): M(u,h)\le t\}
$$
be the growing cluster. For $\WWtil$ we get a 
bound easily by counting self-avoiding paths.

\proclaim{Lemma 5.3} For any finite constant $\beta$, there 
exists a positive $\nu=\nu(\beta)$ such that for any finite
integer $K$, 
$$P\biggl( \text{$\WWtil(t)$ contains a cell $(u,h)$ such
that $|u|+h\ge K$}\biggr) \le 2\exp(-K\beta)
\tag 5.24
$$
as long as $t\le K\nu$. 
\endproclaim

\demo{Proof} Let $\pi$ denote an 
 admissible path fixed to start at the origin $(0,0)=(v^0,k^0)$. 
 The probability in (5.24)
is at most
$$\aligned
&P\biggl( \text{there exists a $\pi$ through at least
$K$ cells with $M(\pi)\le t$} \biggr)\\
\le\;&\sum_{j=K}^\infty
P\biggl( \text{there exists a $\pi$ through exactly 
$j$ cells with $M(\pi)\le t$} \biggr)\\
\le\;&\sum_{j=K}^\infty (2d+1)^j P\bigl( S^1_j\le t\bigr)\\
\le\;&\sum_{j=K}^\infty \exp\lbrakk -j\bigl[
\kappa(t/K)-\log(2d+1)\bigr]\rbrakk\\
\le\;&\sum_{j=K}^\infty \exp(-\beta j) \le 2\exp(-K\beta)\,.
\endaligned
$$
We used again the Cram\'er rate function 
$\kappa(x)=x-1-\log x$ for the Exp($1$) distribution. The 
inequality $P( S^1_j\le t)\le \exp[-j\kappa(t/j)]$ is valid
for $t\le j$. Take $t/K\le\nu$
with $\nu$  small enough so that 
$\kappa(\nu)\ge \beta+\log(2d+1)+1$.
\qed
\enddemo

\proclaim{Lemma 5.4} The processes $\WW(\cdot)$ and $\WWtil(\cdot)$
are equal in distribution.
\endproclaim

\demo{Proof} We can regard both processes as jump processes
on the countable state space of finite connected 
subsets of $\mmZ^d\times\mmZ_+$ that contain the origin. 
Both start from $\{(0,0)\}$. Both processes add new
admissible cells independently at 
rate 1. Comparison of (5.20) and 
(5.23) shows that an admissible new cell, or 
growth cell, is the same for both $\WW(\cdot)$ and $\WWtil(\cdot)$. 
Thus the two processes 
 have identical
infinitesimal rates. \qed
\enddemo

Combining (5.22) and Lemmas 5.3 and 5.4, we get a lower bound for 
$\gamma$: 

\proclaim{Lemma 5.5} There exists a positive constant 
$\nu$ such that $\gamma(x,b)\ge \nu(|x|+b)$ for all
$(x,b)\in\mmR^d\times\mmR_+$. 
\endproclaim

\demo{Proof} Pick $\nu<\nu(1)$ $=$ the constant given by 
Lemma 5.3 for $\beta=1$, and set $t=\nu(|x|+b)$. Then 
by Lemma 5.3 for $K=|[nx]|+[nb]$, 
$$\aligned
P\biggl( R\bigl([nx],[nb]\bigr)\le nt\biggr)
\le\;&P\biggl( Z_{[nx]}(nt) \ge [nb]\biggr)\\
\le\;&P\biggl( W_{[nx]}(nt) \ge [nb]\biggr)\\
\le\;&2\exp\bigl(-|[nx]|-[nb]\bigr).
\endaligned
$$
Thus $R\bigl([nx],[nb]\bigr)\ge nt$ for large enough $n$,
a.s.
\qed
\enddemo

Now we have $\gamma$ bounded both above and below. 
Finally, we convert the limit in Proposition 5.1
to that of Theorem 1. We need one more property for $\gamma$:
 
\proclaim{Lemma 5.6} For $x\in\inte B_0$, there is a unique
finite $h>0$ such that $\gamma(x,h)=1$. 
\endproclaim

\demo{Proof} For this proof we connect 
ballistic deposition on $\mmZ^d\times\mmZ_+$
with first-passage site percolation on $\mmZ^d$. 
Recall the definition of the limit $\mu(x)$ 
in (4.8). By  (4.7) and Proposition 5.1, 
 $\mu(x)=\gamma(x,0)$. By definition,
$B_0=\{x:\mu(x)\le 1\}$, so by homogeneity
$\mu(x)<1$ for $x\in\inte B_0$. Hence for such
$x$ also $\gamma(x,0)<1$. By Lemma 5.5 and the
continuity of $\gamma$ there is some
$h>0$ such that $\gamma(x,h)=1$, so it remains
to rule out the possibility of having $0<h_0<h_1$
such that $\gamma(x,h_0)=\gamma(x,h_1)=1$. But
this and convexity would imply $\gamma(x,0)\ge 1$,
contradicting what was just concluded.
\qed
\enddemo

\demo{Proof of Theorem 1} By the previous lemma, 
a positive function $g$ on $\inte B_0$  is uniquely defined
by the equation $\gamma(x,g(x))=1$. The lower
bound of Lemma 5.5 gives an upper bound for $g$, and $g$ is
concave by the convexity of $\gamma$. A finite
concave function on an open convex set is 
continuous by Theorem 10.1 in \cite{Rf}. By the 
homogeneity of $\gamma$, 
$$\gamma\bigl(x, tg(x/t)\bigr)=t.
\tag 5.25
$$
By the uniqueness in Lemma 5.6 and the monotonicity 
of $\gamma$, $h>tg(x/t)$ [$h<tg(x/t)$] implies
$\gamma(x,h)>t$ [$\gamma(x,h)<t$]. 

Fix $x\in\mmR^d$ and $t>0$ so that $x/t\in\inte B_0$.
Let $\delta_0>0$, and pick $\delta_1\in(0,\delta_0)$. 
Set $b=tg(x/t)-\delta_1$, and pick $\e>0$ so that
$\gamma(x,b)<t-\e$. Fix a large number $m_0$ so that
$m_0(\delta_0-\delta_1)>1$ to take care of the  
effects of 
integer parts. 
$$\aligned
&\lbrakk \liminf_{n\to\infty} n^{-1}Z_{[nx]}(nt) < tg(x/t)-\delta_0
\rbrakk\\
\subseteq\;&\bigcap_{m=m_0}^\infty\bigcup_{n=m}^\infty
\lbrakk Z_{[nx]}(nt) < [ntg(x/t)-n\delta_1]
\rbrakk\\
\subseteq\;&\bigcap_{m=m_0}^\infty\bigcup_{n=m}^\infty
\lbrakk R([nx],[nb]) > n\gamma(x,b)+n\e \rbrakk\\
\subseteq\;
&\lbrakk \limsup_{n\to\infty} n^{-1}
R([nx],[nb]) > \gamma(x,b)+\e \rbrakk\,.
\endaligned
$$
The last event has probability zero by Proposition 5.1. Similarly we
 show that 
$$\limsup_{n\to\infty} n^{-1}Z_{[nx]}(nt) \le
 tg(x/t) \qquad\text{a.s.,}
$$
 and the limit in Theorem 1 is proved. The invariances
of $g$ follow from the corresponding invariances
in the distribution of $Z$. 
\qed
\enddemo

\hbox{}

\head 6. The shape from a translated seed \endhead

This section is a  purely technical extension of
the limit from a seed (Theorem 1) derived in the previous section.
The proof of the hydrodynamic limit (Theorem 2)
 uses the coupling
(3.9) which forces us to consider simultaneously
the whole family $\{Z^v\}$ of processes. Recall
that $Z^v$ stands for the process that grows from a seed
 in cell $(v,0)$, as defined by (3.6). We 
 need a limit where the initial seed and the initial
time point are translated
as the limit is taken. 

In Section 3 we defined  the processes  $\{Z^v\}$ on the 
probability space $(\Omega, \FF,P)$
 of the Poisson jump time processes
$\TT=\{\TT^u\}$. Now we extend  the Poisson processes 
to the entire real line $(-\infty,\infty)$. Since
Poisson points on $(-\infty,0]$ and $(0,\infty)$ are independent, 
this is the same as replacing the original probability
space $(\Omega, \FF,P)$ with a product space 
$$\bigl(\OObar, \FFbar,\Pbar\bigr)=(\Omega^0\times\Omega, \FF^0\otimes\FF,
P^0\otimes P)
$$
where  $(\Omega^0, \FF^0,P^0)$ is the probability space
of the Poisson processes on $(-\infty,0]$. A sample
point of the product space is
$\oobar=(\omega^0,\omega)$, where $\omega$ still
represents  a realization of the i.i.d.\ collection 
$\TT=\{\TT^u:u\in\mmZ^d\}$  of
Poisson point processes on $(0,\infty)$,
while $\omega^0$ represents 
a realization of these processes on  $(-\infty,0]$.
The processes $Z^v$ are defined on 
the space $\OObar$ in the obvious
way, by ignoring the $\omega^0$-component:
$Z^v_u(t)(\omega^0,\omega)=Z^v_u(t)(\omega)$.

Let $\theta_s$ for $-\infty<s<\infty$ denote a 
time translation on the space $\OObar$, so that
the epochs of $\theta_s\oobar$  are those of 
$\oobar$ translated $s$ time units backward.
(In other words, around the time origin 
  $\theta_s\oobar$ looks like $\oobar$ around time point $s$.)
 The random variable $Z^v_u(t)\circ\theta_s$
is computed by reading Poisson jump times from 
time point $s$ onwards and by letting the ballistic deposition
evolve for a duration $t$. The passage times are defined
as before:
$$R^v(u,h)\circ\theta_s=\inf\{t>0: Z^v_u(t)\circ\theta_s\ge h\}\,.
$$
Notice 
 that  $R^v(u,h)\circ\theta_s$ is the {\it amount of time}
it takes to grow up to height $h$, and not the 
time point of the Poisson processes when this happens. 
Note further that if $s<0$ then the random variable 
 $Z^v_u(t)\circ\theta_s$ uses the Poisson processes
also for negative times.

\proclaim{Theorem 4} Let $g$ be the limiting function of
Theorem 1. Then the following holds with probability 1:
For all $x,y\in\mmR^d$, $s\in\mmR$, and $t>0$ such 
that $y\in x+t\cdot \inte B_0$, 
$$\lim_{n\to\infty}\frac1n Z^{[ny]}_{[nx]}(nt)\circ\theta_{ns}
=tg\biggl(\frac{x-y}{t}\biggr)\,.
\tag 6.1
$$
\endproclaim

The remainder of this section proves Theorem 4. 
 We construct step by step an event $\Gamma$ on $\OObar$
that satisfies $\Pbar(\Gamma)=1$
and on which the convergence (6.1) holds.
As in Section 5, we prove the theorem by proving the
convergence of the passage times:
$$\lim_{n\to\infty}\frac{1}{n}R^{[ny]}\bigl([nx],[nb]\bigr)
\circ\theta_{ns}=\gamma(x-y,b)
\tag 6.2
$$
for all $y\in\mmR^d$, $(x,b)\in\mmR^d\times\mmR_+$,
and real $s$.

The reader who is willing to accept the result
may ignore the proof without loss of continuity
and proceed to the next section where the hydrodynamic
limit is proved.

\hbox{}

{\it Step 1.} First we define $\Gamma$ as the event
on which 
$$\lim_{j\to\infty}\frac{1}{2^jm}R^{2^jmv}\bigl(2^jmu,2^jmh\bigr)
\circ\theta_{2^jms_0}=\gamma(u-v,h)
\tag 6.3
$$
for all $v\in\mmZ^d$, $(u,h)\in\mmZ^d\times\mmZ_+$,
rational $s_0$, and all
$m\in\mmN$. 

\hbox{}

We need to argue that $\Pbar(\Gamma)=1$. 
The proof of the Kesten-Hammersley theorem on
p.\ 20--23 of \cite{28} shows that
 convergence along powers of $2$
only 
depends on the distributions of the random variables.
Let $F_n$ denote the distribution of
$R^{nv}(nu,nh)\circ\theta_{ns_0}$.  Then  $F_n$ is also
the distribution of $R(n(u-v),nh)$, and inequality (5.4) shows 
$F_{m+n}\ge F_m\star F_n$. This, and the
existence of second moments, is what is needed 
for the almost sure convergence in (6.3).

\hbox{}

{\it Step 2.} Now define $\Gamma$ to
be the event where the requirement of {\it Step 1 } 
holds, and in addition properties (6.4)--(6.5) below, which
are to hold for all  $y\in\mmQ^d$, $(x,b)\in\mmQ^d\times\mmQ_+$,
rational $s_0$, and rational $s_1,\delta>0$:

(6.4) For large enough $n$, 
$
\bigl|R^{[ny]}\bigl([nx],[nb]\bigr)-R^{[ny]}(v,k)\bigr|\le 2n\delta
$
for all cells $(v,k)$ that satisfy
$\bigl|[nx]-v\bigr|+\bigl|[nb]-k\bigr|\le n\delta.
$

(6.5) For $s_1>0$ there exists $\delta_0(s_1)>0$
 such that if $\delta<\delta_0(s_1)$, 
then for large enough $n$, 
$Z^{[ny]}_v(ns_1)\circ\theta_{ns_0}\ge 0
$
for all sites $v$ such that $|v-[ny]|\le n\delta$.

\hbox{}

Still $\Pbar(\Gamma)=1$. We already argued around (5.15)--(5.17) that (6.4)
can be satisfied almost surely. By Lemma 4.1, 
condition (6.5) can be viewed as 
a percolation question: For each $n$
 define a first-passage
percolation cluster centered at $[ny]$ by 
$$\BB(t)=\{ v\in\mmZ^d: Z^{[ny]}_{v}(t)\circ\theta_{ns_0}\ge 0\}\,.
$$
  If $\delta_0(s_1)$ is chosen so that
 $|v-[ny]|\le n\delta_0(s_1)$ implies $v\in [ny]+n(s_1/2)B_0$,
then (6.5) holds a.s.\ by Corollary 4.2 and the Borel-Cantelli lemma. 

The next step is to improve the convergence in (6.3) to a 
genuine limit on the event $\Gamma$. Fix $m$ for the 
moment. As in (5.8), 
for large enough $n$ there are $j=j(n)$
and $k=k(n)\in\{0,1,\ldots,m-1\}$ such that
$$q(n)\equiv 2^j(m+k)\le n< 2^{j}(m+k+1).
\tag 6.6
$$

Keeping  $v\in\mmZ^d$, $(u,h)\in\mmZ^d\times\mmZ_+$ 
and $s_0\in\mmQ$ fixed, we use
 properties (6.4) and (6.5) of $\Gamma$ to write 
$$\aligned
&R^{q(n)v}\bigl(q(n)u,q(n)h\bigr)\circ\theta_{q(n)(s_0-s)}\\
&\le (n-q(n))s_0+ q(n)s + R^{nv}\bigl(q(n)u,q(n)h\bigr)\circ\theta_{ns_0}\\
&\le (n-q(n))s_0+ns +
 R^{nv}\bigl(nu,nh\bigr)\circ\theta_{ns_0} + 2\delta n\,.
\endaligned
\tag 6.7
$$
This argument will be used several times, so we go over
it once carefully. The first inequality in (6.7) is valid
on the event $\Gamma$ for 
a certain $s=s(m)>0$ that satisfies 
$\lim_{m\to\infty}s(m)=0$,
 and for all large enough $n$, by property (6.5)
for this reason: The ballistic deposition 
process $Z^{q(n)v}\circ\theta_{q(n)(s_0-s)}$ starts 
from a seed in cell $(q(n)v,0)$ when the Poisson 
process clock is at $q(n)(s_0-s)$. 
One way for this 
process  to reach cell $\bigl(q(n)u,q(n)h\bigr)$
is to first spend at most time 
$ (n-q(n))s_0+q(n)s$ to reach cell $(nv,0)$, and from there follow
a new process $Z^{nv}\circ\theta_{ns_0}$ that starts when the Poisson 
process clock is at $ns_0$. Since $0\le n-q(n)\le n/m$,
there is a $\delta=\delta(m)>0$
such that $\lim_{m\to\infty}\delta(m)=0$ and 
$|q(n)v-nv|\le n\delta$ if $n$ is large enough.
Hence by (6.5) we can choose $s=s(m)>0$ such that
  $Z^{q(n)v}\circ\theta_{q(n)(s-s_0)}$ reaches
cell $(nv, 0)$  in time $q(n)s$ if $n$ is large enough.

The second inequality in (6.7) follows from (6.4),
again for  a $\delta=\delta(m)>0$
such that $\lim_{m\to\infty}\delta(m)=0$,  
and for large enough $n$. Also $q(n)\le n$ was used.
A note about terminology: We say that $Z$ reaches 
cell $(u,h)$ in time  $t$ if $Z_u(t)\ge h$. 

Let $n\to\infty$ in (6.7), use the 
limit (6.3) and that
$n-q(n)\le n/m$ to conclude that on $\Gamma$,
$$
\liminf_{n\to\infty}\frac1n R^{nv}(nu,nh)\circ\theta_{ns_0}\ge 
\frac{m-1}m\,\gamma(u-v,h)
-(s+2\delta)-\frac{s_0}m\,.
\tag 6.8
$$
Now let $m\to\infty$, and concurrently we can take $s\to 0$
and $\delta\to 0$. 

To handle the limsup we use the same argument in a reverse
way:
$$\aligned
&R^{nv}\bigl(nu,nh\bigr)\circ\theta_{ns_0}\\
&\le (q(n)-n)s_0 
+q(n)s+R^{q(n)v}\bigl(nu,nh\bigr)\circ\theta_{q(n)(s+s_0)}\\
&\le ns +R^{q(n)v}\bigl(q(n)u,q(n)h\bigr)\circ\theta_{q(n)(s+s_0)}
  + 2\delta q(n)\,.
\endaligned
\tag 6.9
$$
Now the thinking goes like this: 
 The  
process $Z^{nv}\circ\theta_{ns_0}$ starts 
from a seed in cell $(nv,0)$ when the Poisson 
process clock is at $ns_0$. 
One way for this 
process  to reach cell $(nu,nh)$
is to first spend at most time 
$ (q(n)-n)s_0+q(n)s$ to reach cell $(q(n)v,0)$, and from there follow
a new process $Z^{q(n)v}\circ\theta_{q(n)(s+s_0)}$ 
that starts when the Poisson 
process clock is at $q(n)(s+s_0)$. Again, there is a $\delta=\delta(m)>0$
such that $\lim_{m\to\infty}\delta(m)=0$ and 
$|q(n)v-nv|\le n\delta$ if $n$ is large enough.
By (6.5) we can choose $s=s(m)>0$ such that 
$\lim_{m\to\infty}s(m)=0$ and 
  $Z^{nv}\circ\theta_{ns_0}$ reaches
cell $(q(n)v, 0)$  in time $(q(n)-n)s_0+q(n)s$
 if $n$ is large enough.
The second inequality in (6.9) follows from (6.4).

Let $n\to\infty$ in (6.9) to conclude that on the event $\Gamma$,
$$
\limsup_{n\to\infty}\frac1n R^{nv}(nu,nh)\circ\theta_{ns_0}\le 
\frac{m-1}m\,\gamma(u-v,h)
+(s+2\delta)\,.
\tag 6.10
$$
Let $m\to\infty$, and concurrently we can take $s\to 0$
and $\delta\to 0$. 

The limits (6.8) and (6.10) permit us to strengthen the 
definition of $\Gamma$, without losing $\Pbar(\Gamma)=1$:

\hbox{}

{\it Step 3.} The requirements of {\it Step 2}
hold on $\Gamma$, and also 
$$\lim_{n\to\infty}\frac{1}{n}R^{nv}\bigl(nu,nh\bigr)
\circ\theta_{ns_0}=\gamma(u-v,h)
\tag 6.11
$$
for all $v\in\mmZ^d$, $(u,h)\in\mmZ^d\times\mmZ_+$,
and $s_0\in\mmQ$. 

\hbox{}

As previously in Section 5, next we extend the limit 
in (6.11) to all rational sites and cells. Fix
$y\in\mmQ^d$, $(x,b)\in\mmQ^d\times\mmQ_+$,
and $s_0\in\mmQ$. Fix an integer $m$ such that 
 $my\in\mmZ^d$ and  $(mx,mb)\in\mmZ^d\times\mmZ_+$. 
For each $n$, pick $j=j(n)$ so that 
$$q(n)\equiv jm\le n<(j+1)m.$$
Then on $\Gamma$,
$$
\aligned
&\lim_{n\to\infty}\frac{1}{q(n)}R^{q(n)y}\bigl(q(n)x,q(n)b\bigr)
\circ\theta_{q(n)s_0}\\
&=\lim_{j\to\infty}\frac{1}{jm}R^{jmy}\bigl(jmx,jmb\bigr)
\circ\theta_{jms_0}\\
&=\frac1m \gamma(mx-my,mb)\\
&=\gamma(x-y,b)
\endaligned
\tag 6.12
$$
by the homogeneity of $\gamma$. The inequalities (6.7) and (6.9)
can now be repeated, by replacing $v,u,h$ by $y,x,b$, 
and by inserting integer parts where appropriate:
$[ny]$, $[nx]$, $[nb]$. 

We conclude that  $\Pbar(\Gamma)=1$ for the event 
$\Gamma$ defined as follows:

\hbox{}

{\it Step 4.} The requirements of {\it Step 2}
hold on $\Gamma$, and also 
$$\lim_{n\to\infty}\frac{1}{n}R^{[ny]}\bigl([nx],[nb]\bigr)
\circ\theta_{ns_0}=\gamma(x-y,b)
\tag 6.13
$$
for all $y\in\mmQ^d$, $(x,b)\in\mmQ^d\times\mmQ_+$,
and rational $s_0$. 

\hbox{}

The last extension is to 
arbitrary  $y\in\mmR^d$, $(x,b)\in\mmR^d\times\mmR_+$,
and real $s$. Let $\e>0$ be rational, and pick rational 
$\delta\in(0,\e/2)$ so that also $\delta<\delta_0(\e/2)$
[$\delta_0(\e/2)$ as defined in (6.5)].
 Pick rational $y_0,x_0,b_0,s_0$ so that
$$|y_0-y|+|x_0-x|+|b_0-b|+|s_0-s|\le \delta/2\,.$$
Once more we repeat the reasoning that justified
(6.7) and (6.9) to write
$$\aligned
&R^{[ny_0]}\bigl([nx_0],[nb_0]\bigr)\circ\theta_{n(s_0-\e )}\\
&\le  R^{[ny_0]}\bigl([nx],[nb]\bigr)\circ
\theta_{n(s_0-\e ) }+ 2\delta n\\
&\le n(s-s_0+\e ) +
 R^{[ny]}\bigl([nx],[nb] \bigr)\circ\theta_{ns} + 2\delta n\,,
\endaligned
$$
and 
$$\aligned
&R^{[ny]}\bigl([nx],[nb] \bigr)\circ\theta_{ns}\\
&\le  n(s_0-s+\e )
+R^{[ny_0]}\bigl([nx],[nb]\bigr)\circ\theta_{n(s_0+\e )}\\
&\le n (s_0-s+\e )+R^{[ny_0]}\bigl([nx_0],[nb_0]\bigr)
\circ\theta_{n(s_0+\e )}
  + 2\delta n\,.
\endaligned
$$
Let $n\to\infty$ to get 
$$\aligned
&\gamma(x_0-y_0,b_0)-2\delta-(s-s_0+\e )\\
&\le \liminf_{n\to\infty}\frac1n 
 R^{[ny]}([nx],[nb])\circ\theta_{ns}\le \limsup_{n\to\infty}\frac1n 
 R^{[ny]}([nx],[nb])\circ\theta_{ns}\\
&\le \gamma(x_0-y_0,b_0)+2\delta+ (s_0-s+\e ).
\endaligned
$$
As we let $y_0,x_0,b_0,s_0$ approach $y,x,b,s$ we can take
$\delta,\e \to 0$. 
This proves that the limit (6.2) holds on the event $\Gamma$,
and concludes to proof of Theorem 4.

\hbox{}

\head 7. The hydrodynamic limit \endhead

In this section we prove Theorem 2. First we
 work under Assumption A. 
The processes
$\sigma^n$ are constructed according to the 
description of Section 3 on a probability 
space $(\Omega,\FF,P)$  on which are defined 
the Poisson jump time processes $\TT=\{\TT^v\}$
and, statistically independently of $\TT$, 
 the sequence of initial interfaces $\{\sigma^n(0)\}$. 
All processes $\sigma^n$ 
use one and the same  realization
of the Poisson jump time processes $\{\TT^v\}$.
By Corollary 3.1 there
is a single version of the family of ballistic
deposition processes $\{Z^v\}$ grown from seeds as
defined by (3.6) that satisfy
$$\sigma^n_u(t)=\sup_{v\in\mmZ^d}\{ \sigma^n_v(0)
+Z^v_u(t)\}
\tag 7.1
$$
 for all $n\in\mmN$,  $u\in\mmZ^d$ and $t\ge 0$,
 almost surely. Recall that inside the braces the
correct convention is $\infty+(-\infty)=-\infty$. 
 The goal is to prove that 
on some event of full probability, 
 for all $x\in\mmR^d$ and $t>0$, 
$$\lim_{n\to\infty}\frac1n \sigma^n_{[nx]}(nt)=\psi(x,t)\,.
\tag 7.2
$$

For technical reasons we 
extend the Poisson processes to the entire
real line $(-\infty,\infty)$ as was done in
Section 6. Let  $(\Omega^0, \FF^0,P^0)$ be the probability space
of the Poisson processes on $(-\infty,0]$, and consider 
 the product space 
$$\bigl(\OObar, \FFbar,\Pbar\bigr)=(\Omega^0\times\Omega, \FF^0\otimes\FF,
P^0\otimes P)\,.
$$
  A sample
point of the product space $\OObar$ is
$\oobar=(\omega^0,\omega)$, where $\omega\in\Omega$ 
represents  a realization of 
$\bigl(\TT, \{\sigma^n(0)\}\bigr)$ as above, 
while $\omega^0\in\Omega^0$ represents 
a realization of the Poisson  processes on
the negative time line  $(-\infty,0]$.

If we can prove (7.2) in the product space $\OObar$,
it follows for the original space $\Omega$ too. 
For suppose $\Gamma\subseteq\OObar$ is an event such that
$\Pbar(\Gamma)=1$ and (7.2) holds on $\Gamma$. Pick
$\omega^0\in\Omega^0$ such that the $\omega^0$-section 
$$\Gamma_{\omega^0}=\{\omega\in\Omega: (\omega^0,\omega)\in\Gamma\}
$$
satisfies $P(\Gamma_{\omega^0})=1$. Then 
(7.2)  holds on the full-probability
event $\Gamma_{\omega^0}$
because (7.2) depends on $\omega$ only,
and not on $\omega^0$. 
So for the remainder of this section  assume that
all Poisson jump time processes are defined for all 
real times, and we are on the event of full 
probability where Assumption A and  Theorem 4 from
Section 6 
are valid. 

Recall that $\psi$ is defined 
for $(x,t)\in\mmR^d\times(0,\infty)$ by 
$$\psi(x,t)=\sup_{y\in x+ t\cdot\inte B_0}
\lbrak \psi_0(y) +tg\bigl((x-y)/t\bigr)\rbrak\,,
\tag 7.3
$$
and on $\mmR^d\times\{0\}$ set $\psi(x,0)=\psi_0(x)$. 
This first lemma is a consequence of the
 continuity and boundedness of $g$
on $\inte B_0$.

\proclaim{Lemma 7.1} Assume 
 $\psi_0$ is a continuous $[-\infty,+\infty]$-valued
function on $\mmR^d$.
Then  $\psi$ is a continuous $[-\infty,+\infty]$-valued
function on $\mmR^d\times[0,\infty)$.
\endproclaim

\demo{Proof} Let $(x_n,t_n)\to(x,t)$
in  $\mmR^d\times(0,\infty)$. We first 
argue 
$$\liminf_{n\to\infty}\psi(x_n,t_n)\ge \psi(x,t)\,.
\tag 7.4
$$
Suppose $\psi(x,t)>-\infty$ [otherwise (7.4) holds
trivially]. Let $M<\psi(x,t)$, and  
 pick $y\in x+ t\cdot\inte B_0$
so that
$$ \psi_0(y)+tg\bigl((x-y)/t\bigr)>M\,.$$
Since $y\in x_n+ t_n\cdot\inte B_0$ for large enough $n$,
$$\aligned
\liminf_{n\to\infty}\psi(x_n,t_n)
&\ge \liminf_{n\to\infty} \lbrakk 
\psi_0(y) +t_ng\biggl(\frac{x_n-y}{t_n}\biggr)\rbrakk\\
&= \psi_0(y)+tg\bigl((x-y)/t\bigr)\\
&\ge M\,.
\endaligned
$$
The equality above follows from the continuity of $g$
on $\inte B_0$. This proves (7.4).

 To show 
$$M_0\equiv \limsup_{n\to\infty}\psi(x_n,t_n)\le \psi(x,t)\,,
\tag 7.5
$$
pick a subsequence
$n_j$ so that $\lim_{j\to\infty}\psi(x_{n_j},t_{n_j})=
M_0$. 
We may assume $M_0>-\infty$
and $\psi(x,t)<\infty$. Let $M<M_0$,
and for large enough $j$ find  
 $y_j\in x_{n_j}+ t_{n_j}\cdot\inte B_0$ so that 
$$\psi_0(y_j)+t_{n_j}
g\bigl((x_{n_j}-y_j)/t_{n_j}\bigr)>M\,.$$
 Since $B_0$ is compact and  $(x_{n_j},t_{n_j})\to(x,t)$, we may 
pass to a further subsequence so that $y_j\to\ybar$. Define
$$y_j'=x-\frac{t}{t_{n_j}}(x_{n_j}-y_j).$$
Then $y_j'\to \ybar$, $(x_{n_j}-y_j)/t_{n_j}=(x-y_j')/t$, and 
$y_j'\in x+ t\cdot\inte B_0$. 
 We wish to argue that, for large 
enough $j$, 
$$-\infty<\psi_0(y_j)\,,\,\psi_0(y_j')<\infty\,.
\tag 7.6
$$
By continuity of $\psi_0$, this will follow from 
showing that $-\infty<\psi_0(\ybar)<\infty$. We have 
$$\psi_0(\ybar)=\lim_{j\to\infty} \psi_0(y_j')\le \psi(x,t)<\infty\,,
$$
and 
$$\psi_0(\ybar)=\lim_{j\to\infty} \psi_0(y_j)\ge M
-t \|g\|_\infty>-\infty\,.
$$
(7.6) is verified. 
Now for $j$ large enough, 
$$\aligned
M&\le \psi_0(y_j)+t_{n_j}
g\bigl((x_{n_j}-y_j)/t_{n_j}\bigr)\\
&= \psi_0(y_j') +tg\bigl((x-y_j')/t\bigr) +
(t_{n_j}-t) g\bigl((x-y_j')/t\bigr)\\
&\qquad + \psi_0(y_j)- \psi_0(y_j')\\
&\le \psi(x,t) +|t_{n_j}-t|\cdot\|g\|_\infty+ 
\psi_0(y_j)- \psi_0(y_j')\,.
\endaligned
$$
(7.6) was needed to have  the difference $\psi_0(y_j)- \psi_0(y_j')$
 defined and convergent to $0$.
Letting $j\nearrow\infty$ and $M\nearrow 
M_0$ gives (7.5). 
This proves continuity on $\mmR^d\times(0,\infty)$. 
We omit the similar but shorter argument for the 
case  $(x_n,t_n)\to(x,0)$.
\qed
\enddemo

One half of the goal (7.2) follows immediately. 
For arbitrary $y\in x+t\cdot\inte B_0$ such that
$y\in Y_0$, set $u=[nx]$, 
$v=[ny]$ and replace $t$ by $nt$ in (7.1), 
and use assumption (2.5) and Theorem 4 from Section 6  
 to get 
$$\aligned
&\liminf_{n\to\infty}\frac1n \sigma^n_{[nx]}(nt)\\
&\ge \liminf_{n\to\infty}\lbrakk \frac1n \sigma^n_{[ny]}(0)
+ \frac1n Z^{[ny]}_{[nx]}(nt)\rbrakk\\
&=\psi_0(y)+tg\bigl((x-y)/t\bigr)\,.
\endaligned
\tag 7.7
$$
Note that even though the random variable inside the braces
may equal $\infty+(-\infty)$ for finitely many $n$, eventually
$n^{-1} Z^{[ny]}_{[nx]}(nt)$ is finite because its limit 
$tg\bigl((x-y)/t\bigr)$ is finite. 
For each fixed $(x,t)$ take supremum over these admissible $y$'s
to get, by the continuity of the functions involved and by
the denseness of $Y_0$, 
$$\liminf_{n\to\infty}\frac1n \sigma^n_{[nx]}(nt)\ge \psi(x,t).
\tag 7.8
$$
(7.8) holds simultaneously for all $(x,t)$ outside a single event of zero
probability. 

The converse
$$\limsup_{n\to\infty}\frac1n \sigma^n_{[nx]}(nt)\le \psi(x,t)
\tag 7.9
$$
 is where the work is. First we reduce the 
problem to rational $(x,t)$. Suppose (7.9) holds almost surely
 for all
$(x,t)\in\mmQ^d\times\bigl[\mmQ\cap(0,\infty)\bigr]$. 
For rational $(x,t)$ and rational $s>0$, consider the event
$$\aligned
D_{n,s}(x,t)=\lbrak &\text{$\sigma^n_v(nt)\le \sigma^n_{[nx]}(nt+ns)$
for all sites $v$}\\
&\text{such that $|v-[nx]|\le ns/2$} \rbrak\,.
\endaligned
$$
To bound the probability of the complement $D^c_{n,s}(x,t)$, 
fix $n$, $(x,t)$, $s$, and $v$ such that $|v-[nx]|\le ns/2$.
Set $k=\sigma^n_v(nt)$. For times $r\ge nt$,
define 
 a new process $\rho$ by
$\rho(r)=k+Z^v(r-nt)\circ\theta_{nt}$. 
 By the
monotonicity lemma 3.2  
$\rho(nt+ns)\le \sigma^n(nt+ns)$, so 
$\sigma^n_{[nx]}(nt+ns)$ $<$ $k$  implies
that $\rho_{[nx]}(nt+ns)$ $<$ $k$ which in turn 
is equivalent to $ Z^v_{[nx]}(ns)\circ\theta_{nt}$ $=$ $-\infty$.
We can estimate as follows: 
$$\aligned
P\bigl( D^c_{n,s}(x,t)\bigr)
&\le \sum_{v:|v-[nx]|\le ns/2} P\bigl( 
\sigma^n_v(nt)> \sigma^n_{[nx]}(nt+ns) \bigr)\\
&\le \sum_{v:|v-[nx]|\le ns/2} 
P\bigl(Z^v_{[nx]}(ns) =-\infty   \bigr)\\
&\le C_1n^d P\bigl(S^1_{[ns/2]}> ns\bigr)\\
&\le C_1n^d\exp(-C_2n)\,.
\endaligned
$$
We used the fact that the time for process $Z^v $ to
cover cell $([nx],0)$, starting from the seed in cell $(v,0)$, 
is stochastically dominated by $S^1_{|v-[nx]|}$. 
The estimate for 
$P\bigl( D^c_{n,s}(x,t)\bigr) $ 
and Borel-Cantelli imply that, with probability 1,
for any rational $(x,t)$ and $s>0$, the event 
 $D_{n,s}(x,t)$ happens for all large enough $n$. 

Given now arbitrary $(x_0,t_0)$, pick rational 
$(x,t)$ and $s>0$ so that $|x-x_0|<s/4$ and 
$t_0<t$. Then eventually
$\sigma^n_{[nx_0]}(nt)$ $\le$  $\sigma^n_{[nx]}(nt+ns)$, while
by monotonicity $\sigma^n_{[nx_0]}(nt_0)$ $\le$ $\sigma^n_{[nx_0]}(nt)$. 
 We get
$$\limsup_{n\to\infty}\frac1n\sigma^n_{[nx_0]}(nt_0)\le \psi(x,t+s)\,.$$
Let $t+s\searrow t_0$ and $x\to x_0$, use the continuity
Lemma 7.1, and  conclude that now (7.9) holds
for all $(x,t)$ outside a single exceptional event of
zero probability. 

It remains to prove that (7.9) holds a.s.\ for a fixed rational
$(x,t)$. Pick rational $0<s_1<s_0$, with the intention
that  $s_0\searrow 0$ in the end. 
Set $B_1=(t+s_1)\cdot B_0$. Define 
$$\xi_n=\max_{v\in [nx]+nB_1}\{ \sigma^n_v(0)
+Z^v_{[nx]}(nt)\}\,.
\tag 7.10
$$

\proclaim{Lemma 7.2} With probability 1, $\sigma^n_{[nx]}(nt)=\xi_n$
for all large enough $n$.
\endproclaim

\demo{Proof} One way to guarantee the equality 
$\sigma^n_{[nx]}(nt)=\xi_n$ is to require that 
$Z^v_{[nx]}(nt)$ $=-\infty$ for all $v\notin [nx]+nB_1$.
By Lemma 4.1, distributionwise
 this is equivalent to requiring that 
$[nx]\notin \BB^v(nt)$  for all $v\notin [nx]+nB_1$,
where $\BB^v$ is a first-passage percolation cluster
starting from a seed at site $v$, defined as in (2.2)
in terms of $Z^v$. 
Switching to complements, 
$$\aligned
&P\bigl(\sigma^n_{[nx]}(nt)\ne\xi_n\bigr)\\
&\le P\bigl(\text{ 
$[nx]\in \BB^v(nt)$  for some $v\notin [nx]+nB_1$ }\bigr)\\
&= P\bigl(\text{ 
$0\in \BB^v(nt)$  for some $v\notin nB_1$ }\bigr)\\
&\le \sum_{v\notin nB_1} P\bigl(0\in \BB^v(nt)\bigr)\\
&= \sum_{v\notin nB_1} P\bigl( T(0,v)\le nt\bigr)\\
&\le C_1\exp(-C_2n)
\endaligned
$$
for constants $C_i=C_i(t,s_1)$, by Corollary 4.1. The 
conclusion now follows from Borel-Cantelli. 
\qed
\enddemo

 Let
$\{V_i:1\le i\le m\}$ be a collection of closed
neighborhoods whose interiors  cover the compact
set $x+B_1$,
and such  that each 
$$V_i\subseteq y_i+(s_0/2)B_0
\tag 7.11
$$
 for any $y_i\in V_i$, each $V_i$  
lies inside $x+(t+s_0)\cdot\inte B_0$, and satisfies 
assumption (2.6). Since the interiors cover 
$x+B_1$, we have  
$$[nx]+nB_1\subseteq\bigcup_{i=1}^m nV_i
\tag 7.12
$$
for large enough $n$. 
Pick $y_i\in V_i$ such that
$$ g\biggl(\frac{x-y_i}{t+s_0}\biggr)  \le
 \inf_{y\in V_i}g\biggl(\frac{x-y}{t+s_0}\biggr)+s_0.
\tag 7.13
$$ 
By (7.11), Corollary 4.2, and Borel-Cantelli,
 the following holds with 
probability 1: for large enough $n$, 
$$\text{$Z_v^{[ny_i]}(ns_0)\circ\theta_{-ns_0}\ge 0$ 
for all $v\in nV_i$, }
\tag 7.14
$$
 for all $i=1,\ldots,m$. 
In words: Start ballistic deposition processes from 
seeds in cells $([ny_i],0)$ at Poisson process time $-ns_0$.
If $n$ is large enough, 
at Poisson process time $0$ each of these processes  
has grown sufficiently to cover its piece $nV_i$. 
Henceforth assume that we are on this event of full
probability, and that $n$ is large enough for (7.14) to hold.

If we define new processes 
 $\Ztil_u^{[ny_i]}(s)= Z_u^{[ny_i]}(s+ns_0)\circ\theta_{-ns_0}$,
(7.14) gives the inequality
$$\text{$Z_u^v(s)\le \Ztil_u^{[ny_i]}(s)$ 
for all $v\in nV_i$ and $u\in\mmZ^d$}
\tag 7.15
$$
at time $s=0$, for all  $i=1,\ldots,m$. The monotonicity Lemma 3.2
then ensures that (7.15) holds at all times $s\ge 0$, 
and we get 
$$\aligned
\frac1n \xi_n&=\max_{v\in [nx]+nB_1}\lbrakk \frac1n \sigma^n_v(0)
+\frac1n Z^v_{[nx]}(nt)\rbrakk\\
&\le\max_{1\le i\le m}\lbrakk \frac1n \cdot\max_{v\in nV_i}\sigma^n_v(0)
+\frac1n Z^{[ny_i]}_{[nx]}(nt+ns_0)\circ\theta_{-ns_0}\rbrakk\,.
\endaligned
\tag 7.16
$$
Now let $n\to\infty$, apply assumption (2.6), Theorem 4, and 
the choice (7.13) of $y_i$. The limit in Theorem 4 can be taken 
because  $y_i\in V_i\subseteq x+(t+s_0)\cdot\inte B_0$. 
$$\aligned
\limsup_{n\to\infty}\frac1n \xi_n
&\le \max_{1\le i\le m}\lbrakk \sup_{y\in V_i}\psi_0(y)+(t+s_0)
g\biggl(\frac{x-y_i}{t+s_0}\biggr)\rbrakk \\
&\le \sup_{y\in x+(t+s_0)\cdot\inte B_0}\lbrakk \psi_0(y)+(t+s_0)
g\biggl(\frac{x-y}{t+s_0}\biggr)\rbrakk +s_0\\
&= \psi(x,t+s_0)+s_0\,.
\endaligned
\tag 7.17
$$
The argument can be repeated for arbitrarily small $s_0>0$.
Let $s_0\searrow 0$, use the continuity of $\psi$ (Lemma 7.1), 
and then Lemma 7.2 to conclude that (7.9) holds a.s. 
The strong law of Theorem 1 under Assumption A
 is thereby proved. 

Under Assumption B, inequality (7.7), Lemma 7.2, 
and (7.16)--(7.17) 
give, almost surely, 
$$M\le \liminf_{n\to\infty}\frac1n \sigma^n_{[nx]}(nt)
\le \limsup_{n\to\infty}\frac1n \sigma^n_{[nx]}(nt)\le \psi(x,t+s_0)+s_0
$$
for any $M<\psi(x,t)$ and any $s_0>0$. This proves the
statement in Theorem 2 under Assumption B. The argument 
is the same for the weak law under Assumption C. 

\subsubhead Acknowledgements \endsubsubhead 
I thank David Griffeath for an invitation to visit 
University of Wisconsin--Madison where this work was begun, 
and Janko Gravner, Harry Kesten, and Joachim Krug for
valuable comments. 

\hbox{}

\hbox{}

\head References \endhead

\flushpar 
[1] C. Bahadoran (1998). Hydrodynamical limit for spatially
heterogeneous simple exclusion process. \ptrf\ 110 287--331.  

\hbox{}

\flushpar  
[2] I. Benjamini, P. Ferrari, and C. Landim (1996). 
Asymmetric conservative processes with random
rates. \spa\ 61 181--204. 

\hbox{}

\flushpar  
[3] M. G. Crandall, L. C. Evans, and P. L. Lions (1984). 
Some properties of viscosity solutions
of Hamilton-Jacobi equations. Trans. Amer. Math. Soc.
282 487--502.

\hbox{}

\flushpar  
[4] M. G. Crandall and P. L. Lions (1983). Viscosity solutions
of Hamilton-Jacobi equations. Trans. Amer. Math. Soc.
277 1--42.

\hbox{}

\flushpar  
[5] A. De Masi and E. Presutti (1991).  
 Mathematical Methods for Hydrodynamic Limits. 
Lecture Notes in Mathematics 1501,  
 Springer-Verlag,   Berlin.

\hbox{}

\flushpar  
[6] R. Durrett (1995).
Ten  lectures on particle systems. Lecture Notes
in Mathematics 1608 (Saint-Flour, 1993), 97--201. Springer-Verlag. 
   
\hbox{}

\flushpar  
[7] R. Durrett and T. Liggett (1981). The shape of the limit
set in Richardson's growth model. \ap\ 9 186--193.

\hbox{}

\flushpar  
[8] L. C. Evans (1998). Partial Differential Equations. 
American Mathematical Society.

\hbox{}

\flushpar  
[9] D. Griffeath (1979).  Additive
and Cancellative Interacting Particle Systems. 
 Lecture Notes in Mathematics 724, Springer-Verlag. 

\hbox{}

\flushpar  
[10] G. Grimmett and H. Kesten (1984).
First-passage percolation, network flows and
electrical resistances. \zwvg\ 66, 335--366. 

\hbox{}

\flushpar  
[11] T. E. Harris (1972). Nearest-neighbor Markov interaction
processes on multidimensional lattices. Adv.\ Math.\ 9
66--89. 

\hbox{}

\flushpar  
[12] H. Ishii (1984). Uniqueness of unbounded viscosity solutions
of Hamilton-Jacobi equations. Indiana Univ. Math. J. 33 721--748.
  
\hbox{}

\flushpar  
[13] H. Kesten (1986). Aspects of  first-passage
percolation. Lecture Notes in Mathematics 1180,
Springer,  125--264. 

\hbox{}

\flushpar  
[14] H. Kesten (1993). On the speed of convergence in first-passage
percolation. \aap\ 3 296--338.

\hbox{}

\flushpar  
[15] C. Kipnis and C. Landim (1999). Scaling Limit 
of Interacting Particle Systems.
 Grundlehren der mathematischen Wissenschaften, vol 320,
Springer Verlag, Berlin. 

\hbox{}

\flushpar  
[16] J. Krug and P. Meakin (1989). Microstructure and 
surface scaling in ballistic deposition at oblique
incidence. Physical Review A 40 2064--2077.

\hbox{}

\flushpar  
[17] J. Krug and P. Meakin (1991). Columnar growth in  oblique
incidence  ballistic deposition: faceting, noise
reduction, and mean-field theory. Physical Review A 43 900--919.

\hbox{}

\flushpar  
[18] J. Krug and H. Spohn (1991).
  Kinetic
roughening of growing surfaces. 
 Solids far from Equilibrium,
ed. C. Godr\`eche, Cambridge University Press,    
p. 479--582.

\hbox{}

\flushpar  
[19] T. M. Liggett (1985). 
  Interacting Particle Systems.
Springer-Verlag,  New York. 

\hbox{}

\flushpar  
[20] P. Meakin, P. Ramanlal, L. M. Sander, and R. C. Ball
 (1986). Ballistic deposition on surfaces. Physical Review A 34
5091--5103.

\hbox{}

\flushpar  
[21] F. Rezakhanlou (1999). Continuum limit for some
growth models. Preprint.

\hbox{}

\flushpar  
[22] R. T. Rockafellar (1970).  
Convex Analysis. Princeton University Press.          

\hbox{}

\flushpar  
[23] H. Rost (1981). 
Non-equilibrium behaviour of a many particle
process: Density profile and local equilibrium. 
\zwvg\  58 41--53. 

\hbox{}

\flushpar  
[24] T. Sepp\"al\"ainen (1998)
Exact limiting shape for a simplified model of
first-passage percolation on the plane.
 Ann. Probab. 26 1232--1250. 

\hbox{}

\flushpar  
[25] T. Sepp\"al\"ainen (1998). 
Coupling the totally asymmetric simple exclusion process
with a moving interface.
(I Escola Brasileira de Probabilidade, 
IMPA, Rio de Janeiro, 1997), \mprf\ 4 593--628.

\hbox{}

\flushpar  
[26] T. Sepp\"al\"ainen (1999). Existence of hydrodynamics for the 
totally asymmetric simple $K$-exclusion process. 
 \ap\ 27 361--415.

\hbox{}

\flushpar  
[27] T. Sepp\"al\"ainen and J. Krug (1999). Hydrodynamics and 
platoon formation for a  totally asymmetric
exclusion model with particlewise disorder. 
 J.\ Statist.\ Phys. 95  529--571. 

\hbox{}

\flushpar  
[28] R. T. Smythe and J. C. Wierman (1978). First-passage percolation
on the square lattice. Lecture Notes in Mathematics 671, Springer
Verlag. 

\hbox{}

\flushpar  
[29] H. Spohn (1991). 
 Large Scale Dynamics of Interacting Particles.
Springer-Verlag,   Berlin. 

\hbox{}

\flushpar  
[30] M. Talagrand (1995).
Concentration of measure and isoperimetric inequalities in product spaces.
Inst. Hautes \'Etudes Sci. Publ. Math. No. 81 73--205.

\enddocument